\documentclass[reqno]{amsart}%[12pt](tamaño de la letra)
\usepackage{amssymb}
\usepackage{amscd}
\usepackage{amsmath}
\usepackage{enumerate}
\usepackage{graphicx}

\theoremstyle{plain} 
\newtheorem{theo}{Theorem}[section] 
\newtheorem{lem}[theo]{Lemma}
\newtheorem{prop}[theo]{Proposition}
\newtheorem*{nash}{Nash Conjecture}
\theoremstyle{definition} 
\newtheorem{defn}[theo]{Definition} 
\theoremstyle{remark}

\newtheorem{notation}{Notation}

\newenvironment{num}{\medskip
\refstepcounter{subsection}\noindent {\bf (\thesubsection)}}{\vspace{0.2ex}\par} 

\usepackage{hyperref}
\usepackage{comment} 
\includecomment{dibujos}

\usepackage[all]{xy}

\newcommand{\CC}{{\mathbb{C}}}
\newcommand{\DD}{{\mathbb{D}}}

\newcommand{\PP}{{\mathbb{P}}}

\newcommand{\RR}{{\mathbb{R}}}
\newcommand{\SSS}{{\mathbb{S}}}

\newcommand{\ZZ}{{\mathbb{Z}}}

\newcommand{\calA}{{\mathcal{A}}}
\newcommand{\calB}{{\mathcal{B}}}
\newcommand{\calC}{{\mathcal{C}}}
\newcommand{\calD}{{\mathcal{D}}}

\newcommand{\calG}{{\mathcal{G}}}

\newcommand{\calO}{{\mathcal{O}}}

\newcommand{\calR}{{\mathcal{R}}}
\newcommand{\calS}{{\mathcal{S}}}
\newcommand{\calT}{{\mathcal{T}}}
\newcommand{\calU}{{\mathcal{U}}}
\newcommand{\calV}{{\mathcal{V}}}

\newcommand{\calX}{{\mathcal{X}}}

\newcommand{\calZ}{{\mathcal{Z}}}
\newcommand{\G}{\mathrm B\mathbf{I}}
\newcommand{\I}{\mathbf{I}}
\newcommand{\X}{\CC^2/\mathrm B\mathbf{I}}
\newcommand{\ori}{\{O\}}

\newcommand{\comp}{{\circ}}
\newcommand{\cont}{{\subset}}

\newcommand{\mm}{\mathfrak{m}}

\newcommand{\codim}{{\mathrm{codim}}}
\newcommand{\mult}{{\mathrm{mult}}}

\newcommand{\inv}{^{-1}}

\begin{document}
\title{Nash Problem for Quotient Surface Singularities.}
\author{Mar\'ia Pe Pereira}
\address{Departamento de \'Algebra\\Universidad Complutense de Madrid\\
Parque de Ciencies 3\\
28040 Madrid-Espa\~na}
\curraddr{Institut de Math\'ematiques de Jussieu\\
\'Equipe G\'eom\'etrie et Dynamique\\
175 rue du Chevaleret\\
75013 Paris} 
\email{maria.pe@mat.ucm.es}
\thanks{The author is supported by Caja Madrid. In the course of this research the author was supported by the ERC Starting Grant project TGASS in ICMAT-UCM-UAM-Carlos III and by a FPI grant associated to Spanish Project MTM2004-08080-C02-01 in Universidad Complutense de Madrid.}
\date{15-11-2010}
%{Primary: 32, 32}

\begin{abstract}
We give an affirmative answer to Nash problem for quotient surface singularities.
\end{abstract}

\maketitle

%%%%%%%%%%%%%%%%%%%%%%%%%%%%%%%%%%%%%%%%%%%%%%
\section{Introduction}\label{sec:intro}
%%%%%%%%%%%%%%%%%%%%%%%%%%%%%%%%%%%%%%%%%%%%%%
In 1968 (published as \cite{Nash}), J. Nash introduced the study of arc spaces of an algebraic or analytic variety. Arc spaces in a variety $X$ have an infinite dimensional algebraic variety structure, viewed as a limit of spaces of truncated arcs, which have a natural structure of finite-dimensional algebraic variety inherited from $X$. Nash proved in \cite{Nash} that arc spaces have finitely many irreducible components.

In \cite{Nash}, he forecasts a relation between the irreducible components of the arc space $\calX_\infty$ of a variety $X$ and the essential divisors  of a resolution of  singularities $\pi:\tilde{X}\to X$ (the divisors that appear in any resolution up to birational maps).
In fact he defined a natural mapping from the set of irreducible components of the arc space to the set of essential divisors. He proposed its study, proved that it is injective and conjectured that it is a bijection in the surface case. For a more extended introduction to the problem one can see  \cite{Nash} or \cite{IK}.

In dimension $4$, Ishii and Kollar in \cite{IK} found an example of singularity with non-bijective Nash mapping. On the other hand the bijectivity
of the Nash mapping is still open for surfaces and has been proved in many classes of singularities: $\mathbf{A}_k$ singularities \cite{Nash}; normal minimal surface singularities \cite{Reg-min} (and also \cite{Cam-min} and \cite{Fer}), sandwiched singularities \cite{Reg-sand}, \cite{MonReg-sand}; the dihedral singularities $\mathbf{D}_n$ \cite{Cam-Dn}; germs with a good $\CC^*$-action such that $Proj\calS$ is not rational, following \cite{Mon-LN} and \cite{Reg-sand}; a family of non-rational surface singularities, given by the property $E\cdot E_k<0$ for any $i$ in \cite{CamPPP1}; toric singularities \cite{IK}; quasi-ordinary singularities \cite{Ishii-qo}, \cite{Pedro}; stable toric varieties \cite{Petrov}; some other families of higher dimensional singularities \cite{CamPPP2}; ...

Besides, there are some papers proving reductions of the problem. For instance, in \cite{MonReg-uniruled} it is proved that divisors that
are not uniruled are in the image of the Nash map and that the so-called \emph{lifting wedge property} for the case of quasirational surfaces would imply the surjectivity of the Nash map for normal surface singularities. In \cite{Bob} it is proved that Nash problem for surfaces only depends on the topological type of the surface singularity and that a positive answer for the case of rational homology sphere links implies a positive answer for any normal surface;  
 %and showing that certain kind of components must be in the image of Nash mapping, no
In \cite{Reg-CS}, \cite{MonReg-uniruled} and \cite{Bob} certain characterizations of the bijectivity of Nash mapping are proved in terms of uniparametric families of arcs.% (we will explain more below on this approach).%, in \cite{MonReg-uniruled} it is proved that divisors that
%are not uniruled are at the image of Nash mapping, in \cite{Bob} it is proved that Nash problem for surfaces is a topological problem, and that a positive answer for the case of rational homology sphere links implies a positive answer for any normal surface.

Despite these advances, up to this moment, Nash Problem %has proved to be very difficult and 
is still open for the very basic cases of surface rational double points $\textbf{E}_6$, $\textbf{E}_7$ and $\textbf{E}_8$ (at the time of submitting this paper a proof for $\textbf{E}_6$ appeared in \cite{CamSpi}). In this paper we prove the following
\begin{theo}\label{theo:nash_quot}
Nash mapping is bijective for any quotient surface singularity, in particular for $\textbf{E}_6$, $\textbf{E}_7$ and $\textbf{E}_8$.
\end{theo}    

Our proof is based on the characterization of bijectivity given in \cite{Bob} in terms of families of convergent arcs depending holomorphically on a parameter. The convergence allows us to work geometrically and to use topological arguments which are essential reduction steps for the proof.  

Let us sketch briefly the strategy of the proof.
In the surface case, there exists a minimal resolution $\pi:(\tilde{X},E)\to (X,O)$, where the exceptional divisor $E$ decomposes into a finite number of irreducible components $E=\cup_{k=1}^r E_k$, which are the essential components.
Given a divisor $E_k\cont\tilde{X}$ one can consider the set $N_k$ of arcs whose lifting to the resolution meets $E_k$. Nash proved in \cite{Nash} the irreducibility of the closure of the sets $N_k$; hence, the bijectivity of Nash mapping is equivalent to the fact that no $N_i$ is in the Zariski closure of a different $N_j$. If an inclusion $N_i\subseteq \overline{N}_j$ for $i\neq j$ occurs, we say that there is an \emph{adjacency} (from $N_j$ to $N_i$). 

Using the shape of the resolution graphs of the quotient surface singularities given in \cite{Bri} and a comparison theorem of \cite{Bob}, we can reduce the proof to the singularities 
$\textbf{D}_n$ and $\textbf{E}_8$ (in fact a version of a result in \cite{Cam-min} is enough for this reduction). The singularities $\textbf{D}_n$ were settled in \cite{Cam-Dn} and we also recover its result using our method in \cite{tesis}. % but we remark that they follow easily using our method. 
We apply several criteria
to rule out adjacencies at different levels of depth.
In fact, suffice already the simplest criteria here  to prove the conjecture for most of quotient surface singularities. For the sake of brevity we include here the complete proof only for the $\textbf{E}_8$ singularity. We briefly indicate in Section \ref{sec:quot} which of the criteria allow to prove Nash conjecture for quotient surface singularities without using the comparison theorem in \cite{Bob}. The more elaborate criteria are needed only for the $\textbf{E}_k$'s. A detailed exposition of this method for general quotient surface singularities can be found in \cite{tesis}.  

The idea is to use families of arcs and to pull them back by the quotient mapping to families of curves in $\CC^2$. At the first level of depth we use upper semicontinuity of intersection multiplicity in $\CC^2$. This criterion seems to be equivalent to the \emph{valuative criterion} in \cite{Reg-min} and \cite{Cam-min} that compares orders of functions on the divisors $E_k$. % this rules out the same adjacencies than PLENAT rules out in. 
At a second level, using topology of plane curves singularities, we discover a new phenomenon, the so called {\em returns}, not used before in Nash Problem, but which is the 
key for our proof of the impossibility of the most difficult adjacencies. A \emph{return} occurs in a family of arcs when the generic members of the family go through the origin of the surface singularity more than once. Combining returns with upper semicontinuity of intersection multiplicity we rule out a second set of adjecencies. For the last adjacencies we find implicit equations for the families of arcs, use versal deformations, and find a surprising proof in terms of dimensions of suitable strata in the base of the versal deformation.

%%%%%%%%%%%%%%%%%%%%%%%%%%%%%%%%%%%%%%%%%%%%%%
\section{Preliminaries}\label{sec:pre}
%%%%%%%%%%%%%%%%%%%%%%%%%%%%%%%%%%%%%%%%%%%%%%
Let $(X,O)$ be a normal surface singularity.
\begin{defn}
An \emph{arc} is a morphism 
$$Spec(\CC[[t]])\to X$$
sending the special point to the origin $O$. We denote by $\calX_\infty$ the space of arcs of $X$. 
An arc is said to be \emph{convergent} if it is given by a holomorphic mapping 
$$\gamma:(\CC,0)\to X.$$%the power series defining it are convergent.
\end{defn}

Let 
$$\pi:(\tilde{X},E)\to (X,O)$$
be the minimal resolution of singularities, with exceptional divisor $E$. Let $E=\cup_{k=1}^rE_k$ be the decomposition of $E$ in irreducible components.
By the valuative criterion of properness any arc $\gamma$ admits a unique lifting
$$\tilde{\gamma}:Spec(\CC[[t]])\to\tilde{X}.$$
Of course, if $\gamma$ is convergent, so is $\tilde{\gamma}$. 
\begin{defn}
An arc is said \emph{to have a transverse lifting} if its lifting $\tilde{\gamma}$ meets only 
one irreducible component $E_k$ transversely at a smooth point of $E$.
\end{defn}

Given a divisor $E_k$, define
$$N_k=\{\gamma\in \calX_\infty: \tilde{\gamma}(0)\in E_k\},$$
$$\dot{N}_k=\{\gamma\in N_k: \tilde{\gamma}\rm{\ is \, transverse\ to}\, E_k \}.$$
The set $\dot{N}_k$ is dense in $N_k$ and the closure set $\overline{N}_k$ is irreducible (\cite{Nash}). On the other hand the arc space $\calX_\infty$ is the union of these sets. Hence, the Nash problem asks whether all of these sets are in fact irreducible components of the space of arcs or  there is any inclusion $\overline{N}_i\subseteq\overline{N}_j$ for $i\neq j$. Such an inclusion $\overline{N}_i\subseteq\overline{N}_j$ is called an \emph{adjacency}. 
\begin{nash} Given a normal surface singularity, for any two essential divisors $E_i$ and $E_j$ we have 
$$\overline{N}_i\nsubseteq\overline{N}_j.$$
\end{nash}
We will study Nash problem via a characterization using wedges, which are uni-parametric
families of arcs.
\begin{defn}
A \emph{wedge} is a morphism
$$\alpha:Spec[[t,s]]\to X$$
such that the closed subset $V(t)$ is sent to the origin $O$. A wedge is said to be \emph{convergent} if it is given by a holomorphic mapping 
$$\alpha:(\CC^2,O)\to X.$$
We denote $\alpha_0:=\alpha|_{(\CC\times\{0\},0)}$ and $\alpha_s:=\alpha|_{(\CC\times\{s\},0)}$ for any $s\in (\CC,0)\setminus\{0\}$. Unless we state the contrary, $s$ stands for any small enough parameter value different from $0$.% for any value of the parameter $s\in \CC$. 
\end{defn}

Our proof uses the following characterization:

\begin{theo}\label{theo:ja}{\cite{Bob}} Let $E_i$ and $E_j$ be two components of $E$. The following are equivalent:
\emph{\begin{enumerate}
	\item $N_i\cont\overline{N}_j$;
	\item there exists a convergent wedge $\alpha$ 
	with $\alpha_0\in \dot{N}_i$ and $\alpha_s\in \dot{N}_j$;
	%whose special arc $\alpha_0$ has lifting to $\tilde{X}$ transverse through $E_u$ and whose generic arc has lifting transverse through  $E_j$;
	\item for any convergent arc $\gamma$ in $\dot{N}_i$
	%whose lifting to $\tilde{X}$ is transverse through the divisor $E_u$, 
	there exists a convergent wedge $\alpha$ with $\alpha_0=\gamma$ and $\alpha_s\in\dot{N}_j$.%generic arc with lifting transverse through $E_j$.
\end{enumerate}}
\end{theo}

The following result reduces Theorem \ref{theo:nash_quot} to the cases of $\textbf{E}_8$ (and $\textbf{D}_n$):
\begin{theo}[\cite{Bob}]
\label{theo:comparacion1}
The following statements hold:
\begin{itemize}
\item Let $\calG_1$ be a graph contained in a negative definite weighted graph $\calG_2$.
If Nash mapping is bijective for a surface singularity with resolution graph $\calG_2$, then it is bijective also for a surface singularity with resolution graph $\calG_1$.
\item Let $\calG_1$ be a graph obtained from a negative definite weighted $\calG_2$ by decreasing the self-intersection weights of vertices (the graph $\calG_1$ is automatically negative definite).
If Nash mapping is bijective for a surface singularity with resolution graph $\calG_2$, then it is bijective also for a surface singularity with resolution graph $\calG_1$.
\end{itemize}
\end{theo}
In fact, a version of a result in \cite{Cam-min} is enough to make the reduction in our case.

%%%%%%%%%%%%%%%%%%%%%%%%%%%%%%%%%%%%%%%%%%%%%%%%%%%%%
%%%%%%%%%%%%%%%%%%%%%%%%%%%%%%%%%%%%%%%%%%%%%%%%%%%%%
\section{The Icosahedral Singularity $\textbf{E}_8$.}\label{sec:E_8}
%%%%%%%%%%%%%%%%%%%%%%%%%%%%%%%%%%%%%%%%%%%%%%%%%%%%%
%%%%%%%%%%%%%%%%%%%%%%%%%%%%%%%%%%%%%%%%%%%%%%%%%%%%%
%\subsection{
%%%%%%%%%%%%%%%%%%%%%%%%%%%%%%%%%%%%%%%%%%%%%%%%%%%%%
\begin{num}\textbf{ Definition and some properties.} \label{def:E_8}
%%%%%%%%%%%%%%%%%%%%%%%%%%%%%%%%%%%%%%%%%%%%%%%%%%%%%
Let $\mathbf{I}$ be the group of isometries of the icosahedron that preserve the orientation. It is a finite subgroup of $SO(3,\RR)$ and then it leaves $\SSS^2$ invariant. By stereographic projection we identify $\PP^1$ with $\SSS^2$ and hence $SO(3,\RR)$ with $PSL(2,\CC)$. We see $\mathbf{I}$ as a subgroup of $PSL(2,\CC)$.  
The \emph{binary icosahedral group} $\mathrm B\mathbf{I}$ is the subgroup of $SL(2,\CC)$ which is the preimage of the icosahedral group $\I$ by the quotient map:
$$\rho: SL(2,\CC)\to PSL(2,\CC)=SL(2,\CC)/\{Id,-Id\}.$$
The group $\mathrm B\mathbf{I}$ acts properly discontinuously in $\CC^2$ and then, according to \cite{Car}, the quotient surface $\CC^2/\mathrm B\mathbf{I}$ is a normal surface, called the $\mathbf{E_8}$ singularity. Since no element leaves a hyperplane fixed, the projection map 
$$p:\CC^2\to \CC^2/\mathrm B\mathbf{I}$$ is a regular covering map outside the origin. All over the paper we denote $\X$ by $X$. 

Given a group $G$ acting in $\CC^2$, an equation $F\in \CC\{u,v\}\setminus\ori$ is called $G$-\emph{semi-invariant} if for each $g\in \mathrm{G}$ there exists $\chi_F(g)\in \CC^*$ such that 
\begin{equation}\label{eq:G-inv}
F\comp g=\chi_F(g)F.
\end{equation}
The mapping $\chi_F:\mathrm{G}\to \CC^*$ is a homomorphism and is called \emph{the character of F}. By \cite{Bri} the ring of functions of $\textbf{E}_8$ is a Unique Factorization Domain. This implies (by \cite{Mat}, pag. 141) that the character of any $\mathrm B\mathbf{I}$-semi-invariant function in $\CC[u,v]$ is $1$, that is, it is in fact  $\mathrm B\mathbf{I}$-invariant. That the only possible character is $1$ is also obtained in a direct way in \cite{Lam}. 

The action of $\mathrm{G}$ on the $\CC$-algebra $\CC\{u,v\}$ given by (\ref{eq:G-inv}) is linear and hence it respects the decomposition of $F$ into its homogeneous terms. Particularly $F$ is G-invariant if and only if its homogeneous terms are G-invariant. By the Fundamental Theorem of Algebra every form in two variables can be written as
\begin{equation}\label{eq:form}
H(u,v)=(b_1u-a_1v)^{r_1}...(b_mu-a_mv)^{r_m}
\end{equation}
where $[a_1:b_1]$,...,$[a_m:b_m]\in \PP^1$. 

The group $\mathrm B\mathbf{I}$ acts in such a form as the group $\mathbf{I}$ acts in the points $[a_1:b_1]$,...,$[a_m:b_m]\in \PP^1$ (we identify $\PP^1$ with the unit sphere $\SSS^2$). The form $H$ is $\mathrm B\mathbf{I}$-invariant if and only if this set of points is invariant by $\mathbf{I}$.
The elementary $\mathbf{I}$-orbits are the orbits of  regular points $q$ -which have 60 points and are called $P_q$- and the exceptional orbits -which correspond to the 30 mid-edges, 20 face-centers and 12 vertices of the icosahedron-. They are denoted by $E$, $F$ and $V$ and have degrees 30, 20 and 12 respectively. 

It can be checked that $E$, $F$ and $V$ form a minimal set of generators of the $\mathrm B\mathbf{I}$-invariant forms and that they give an expression for the quotient map 
$$p:\CC^2\longrightarrow \CC^3$$
\begin{equation}\label{eq:p}(u,v)\mapsto (E,F,V)
\end{equation}
whose image has equation $x^2+y^3+z^5=0$. The details can be found in \cite{Lam}.
\end{num}

%%%%%%%%%%%%%%%%%%%%%%%%%%%%%%%%%%%%%%%%%%%%%%%%%%%%%
%\subsection{Minimal resolution.}
\begin{num}\textbf{ Minimal resolution.}\label{sec:min_res}
%%%%%%%%%%%%%%%%%%%%%%%%%%%%%%%%%%%%%%%%%%%%%%%%%%%%%
 We construct following \cite{Bri} the minimal resolution of $\X$. 
Let us take the blow-up of the origin of $\CC^2$
$$\pi_0:Bl_0\CC^2 =\{(x,y;[a:b])\in\CC^2\times \PP^1: xb=ya \}\to \CC^2$$
$$(x,y;[a:b])\mapsto (x,y).$$
Let us denote the exceptional divisor by
$$F_1:=\{(0,0;[a:b])\in\CC^2\times \PP^1\}.$$
The action of $\G$ over $\CC^2$ admits a lifting to $Bl_0\CC^2$ assigning to any $g\in \G$ the mapping, denoted by abuse of notation also by $g$, the following mapping
\begin{equation}
g:Bl_0\CC^2\to Bl_0\CC^2, \ g((x,y;[a:b]))=(g(x,y);[g(a,b)]).
\label{eq:lifted_action}
\end{equation}
In particular, the lifted action leaves $F_1$ invariant and acts on it as the group $\I$. We denote the quotient map by
\begin{equation}\label{eq:p_1}p_1:Bl_O\CC^2\to Bl_O\CC^2/\G\end{equation}
and define
\begin{equation}\label{eq:E_1}E_1:=p_1(F_1).\end{equation}

The quotient space $Bl_0\CC^2/\G$ is smooth except at 3 points which  are the image of the points with stabiliser different from $\ZZ_2$. These points correspond to the 30  mid-edges, 20 face-centers and 12 vertices of the icosahedral and have stabilisers $\textbf{C}_{4,2}$, $\textbf{C}_{6,4}$ and $\textbf{C}_{10,8}$ respectively where $\textbf{C}_{n,q}$, with $q<n$, denotes the cyclic group action in $\CC^2$ given by \begin{equation}\label{eq:Xnq}\left(\begin{array}{cc}\zeta_{\tiny{n}}&0\\0&\zeta_{\tiny{n}}^q\end{array}\right)\end{equation}
with $\zeta_{\tiny{n}}$ an $n$-th root of unity. The images of the special orbits are the  singularities $\textbf{A}_1$, $\textbf{A}_2$ and $\textbf{A}_4$ respectively. The projection $Bl_0\CC^2\to \CC^2$ is $\G$-equivariant and induces the map
 $$q:Bl_0\CC^2/\G\to \CC^2/\G=X$$
which is a proper holomorphic map that is biholomorphic out of the origin of $X$. Then, resolving the 3 singularities of $Bl_0\CC^2/\G$ we get a resolution of $\X$ that can be checked that it is the minimal one, with resolution graph as in Figure \ref{fig:graph}, and that we denote as follows:
$$\pi:\tilde{X}\to X=\X.$$

\begin{dibujos}
\begin{figure}[!h]
\centering
\includegraphics[width=2.6in]{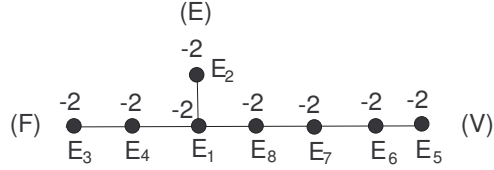}
\caption{Resolution graph of $\textbf{E}_8$.}\label{fig:graph}
\end{figure}
\end{dibujos}
\end{num}
%%%%%%%%%%%%%%%%%%%%%%%%%%%%%%%%%%
%\subsection{
\begin{num}\textbf{ Resolution of indeterminacy.}\label{sec:indet} 
%%%%%%%%%%%%%%%%%%%%%%%%%%%%%%%%%%%%%%%%%%%%%%%%%%%%%
We will need the following construction. Given the mappings $p$ and $\pi$ as above, by \emph{resolving the indeterminacy} of $\pi\inv\comp p:\CC^2--->\tilde{X}$ we mean finding a smooth variety $Y$, a proper birational morphism $\bar{\pi}:Y\to \CC^2$ and a holomorphic morphism $\bar{p}:Y\to \tilde{X}$ that close the following diagram
\begin{equation}\label{dia:indet}\xymatrix{
Y\ar_{\bar{\pi}}[dd]\ar^{\bar {p}}[rr]&&\tilde{X}\ar_{\pi}[dd]\\
&&\\
\CC^2\ar@{.>}[uurr] \ar^{p}[rr]&&X.
}.
\end{equation}
We obtain the construction inductively, gluing the resoltuions of indeterminacy for the singularities $\textbf{A}_1$, $\textbf{A}_2$ and $\textbf{A}_4$. We denote them as for the case of $\X$ but with a subindex indicating the action of the stabiliser.  
For instance, for $\textbf{A}_2$, with stabiliser $\textbf{C}_{6,4}$, we consider the quotient map 
\begin{eqnarray}\label{eq:p_6,4} p_{6,4}: \CC^2&\to& \CC^3\nonumber\\
(u,v)&\mapsto &(u^6,v^3,u^2v).
\end{eqnarray}
The image is denoted by $X_{6,4}:=\CC^2/\textbf{C}^{6,4}$ and has equation $xy-z^3=0$. Performing one blow up in the origin of $\CC^3$ we get the resolution of singularities $\pi_{6,4}:\tilde{X}_{6,4}\to X_{6,4}$ inside $\CC^3\times \PP^2$. The exceptional divisor is the union of 2 transverse irreducible curves of genus $0$, both with self-intersection $-2$. In order to get a resolution of indeterminacy, we make a sufficient sequence of blow-ups in points of the source $\CC^2$ to kill the indeterminacy from these spaces to $\tilde{X}_{6,4}\cont\CC^3\times \PP^2$ as in the diagram (\ref{dia:indet}). We finally obtain a space $Y_{6,4}$ and two mappings: the composition of the blow-ups  $\bar{\pi}_{6,4}:Y_{6,4}\to \CC^2$ and the resolution of indeterminacy $\bar{p}_{6,4}:Y_{6,4}\to \tilde{X}_{6,4}$. In particular, for any component $E_k$ of the exceptional divisor of $\tilde{X}_{6,4}$ there is an irreducible component $\bar{E}_k$ of the divisor $\bar{\pi}_{6,4}\inv(O)$ that is mapped onto it. These divisors are called \emph{dicritical}.
The restriction of $\bar{p}_{6,4}$ around a smooth point of $\bar{E}_k$ is of the form $(u,v)\mapsto (u^m,v)$ where $\bar{E}_k$ is given by the equation $u=0$. Then, we get that the restriction $\bar{p}_{6,4}|_{\bar{E}_k}$ has degree $|\textbf{C}_{6,4}|/m=6/m$. In Figures \ref{fig:A1_indet}-\ref{fig:A4_indet} we draw the resolutions and the resolutions of indeterminacy for $\textbf{A}_1$, $\textbf{A}_2$ and $\textbf{A}_4$: observe that we denote the essential divisors for the singularities $\textbf{A}_i$ by the name that they have when we glue their resolutions to construct the resolution of $\textbf{E}_8$. We denote by  $\bar{E}_k$ the corresponding dicritical divisors.

We also find a regular parametrization of generic arcs in $\CC^2$ whose lifting to $\tilde{X}$ are transverse to the dicritical divisors. For example in $A_2$, we can express locally the mapping $\bar{\pi}_{6,4}$ around a generic point of $\bar{E}_4$ as $(u,v)\mapsto (uv,v)$ and around a generic point of $\bar{E}_3$ as $(u,v)\mapsto (u,u^4v)$. Then, the image in $\CC^2$ of an arc $(\lambda,t)$ in the chosen chart around $\bar{E}_4$ is $(\lambda t,t)$ and the image of $(t,\lambda)$ in the chart around $\bar{E}_3$ is $(t,\lambda t^4)$. 

In order to get a resolution of indeterminacy $Y$ for $X$ as in the diagram \ref{dia:indet}, we glue the resolution of indeterminacy $Y_{4,2}$, $Y_{6,4}$ and $Y_{10,8}$ obtained for the cyclic singularities around their preimages in $F_1$. We have to glue 30 copies of the resolution of indeterminacy of $A_1$, 20 of the one of $A_3$ and 12 of the one of $A_4$. Consequently, a divisor $E_k$ of $\tilde{X}$ will have a different number of dicritical divisors over it, 30, 12 or 20, depending on the cyclic singularity which it comes from. We denote the dicritical divisors over $E_k$ by $\bar{E}_{k}^1$, ..., $\bar{E}_k^n$. Given a divisor $\bar{E}_k^l$, the degree of the covering maps from it to $E_k$ is the same than in the corresponding cyclic case and consequently coincides for all the dicritical divisors over $E_k$.  
\begin{dibujos}
\begin{figure}[!h]
\centering
\includegraphics[width=1.7in]{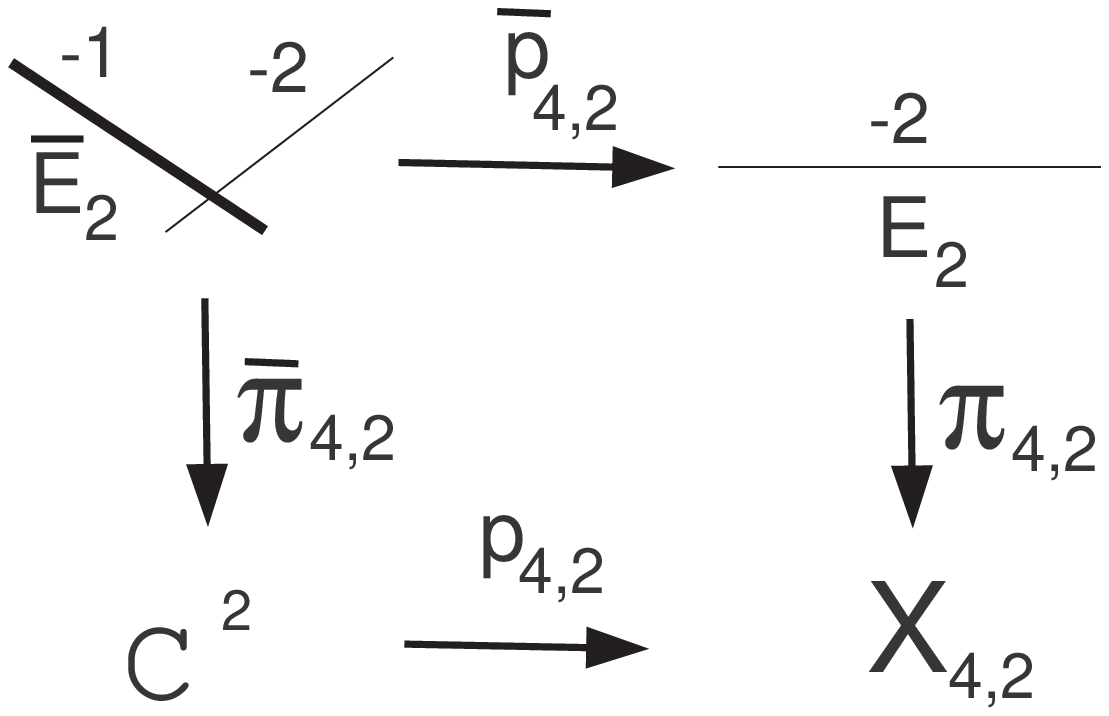}
\rule[0mm]{18mm}{0mm}%{\rule{30pt}{1cm}}
\includegraphics[width=2.6in]{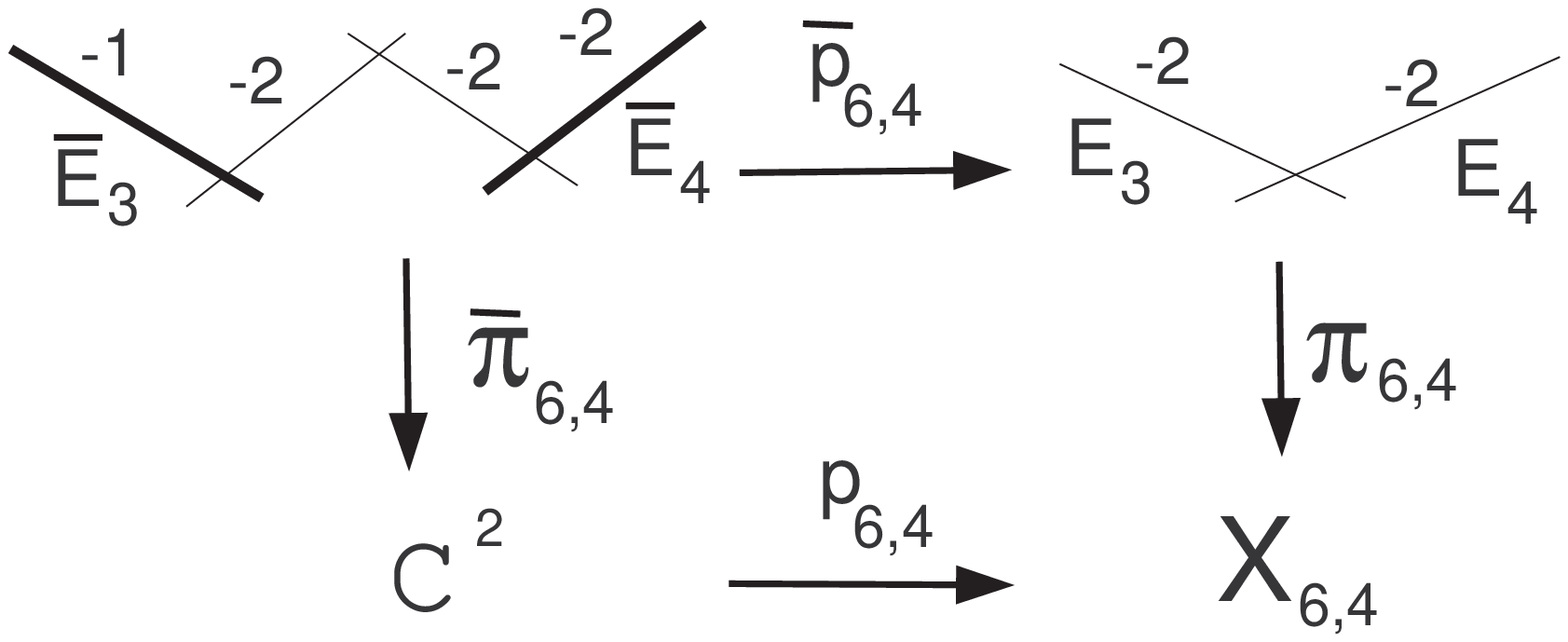}
\caption{Resolution of indeterminacy of $\textbf{A}_1$ and $\textbf{A}_2$.}\label{fig:A1_indet}
\end{figure}
\end{dibujos}

\begin{dibujos}
\begin{figure}[!h]
\centering
\includegraphics[width=5in]{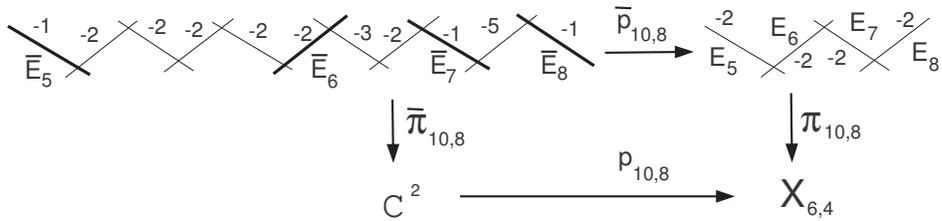}
\caption{Resolution of indeterminacy of $\textbf{A}_4$.}\label{fig:A4_indet}
\end{figure}
\end{dibujos}

The image in $\CC^2$ of an arc in $Y$ transverse to a dicritical divisor over some $E_k$, for $k\neq 1$, can be obtained from the ones obtained in the cyclic case. The parametrizations obtained for the cyclic case can be considered to live in $Bl_O\CC^2$. Then, we do one more blow down to get to $\CC^2$. The tangent line to such an arc is given by the point in $F_1\cont Bl_O\CC^2$ that the parametrization (of its lifting) in $Bl_o\CC^2$ meets. 

In Table \ref{tabla1} we give a summary of this information for each divisor of $\tilde{X}$.

\begin{table}[!h]%
\centering
\caption{Some information about the resolution of indeterminacy $Y$. } 
\label{tabla1}
\begin{tabular}{|l|l|l|l|l|l|l|l|l|}
\hline
Divisor&$E_1$&$E_2$&$E_3$&$E_4$&$E_5$&$E_6$&$E_7$&$E_8$\rule[-2.25mm]{0mm}{7mm}\\
\hline 
Stabiliser in $Bl_O\CC^2$&$\textbf{C}_{4,2}$&$\textbf{C}_{6,4}$&$\textbf{C}_{6,4}$&$\textbf{C}_{10,8}$&$\textbf{C}_{10,8}$&$\textbf{C}_{10,8}$&$\textbf{C}_{10,8}$&$\textbf{C}_{10,8}$\rule[0mm]{0mm}{5mm}\\
Orbit in $F_1$&$P_q$&E&F&F&V&V&V&V\\
Dicritical divisors&1&30&20&20&12&12&12&12\\
Degree of $\bar{E}^l_k\to E_k$&60&1&1&2&1&2&1&2\\
Parametrizations &$(t,\lambda t)$&$(t,\lambda t^3)$&$(t,\lambda t^5)$&$(t,\lambda t^2)$&$(t,\lambda t^9)$&$(t,\lambda t^4)$&$(t^3,\lambda t^7)$&$(t^2,\lambda t^3)$\rule[-2.25mm]{0mm}{0mm}\\
\hline
\end{tabular}
\end{table}
\end{num}
%%%%%%%%%%%%%%%%%%%%%%%%%%%%%%%%%%%%%%%%%%%%%%%%%%%%%
\section{Proof of Nash Conjecture for $\textbf{E}_8$.}
%\subsection{}
%%%%%%%%%%%%%%%%%%%%%%%%%%%%%%%%%%%%%%%%%%%%%%%%%%%%%
The proof consists in finding a contradiction to the existence of each one of the 56 possible adjacencies $N_i\subseteq\overline{N}_j$. Thanks to Theorem \ref{theo:ja}, we can work only with convergent arcs and wedges. Thus, all the arcs and wedges that appear in the proof are convergent.

\subsection*{Part I}

%%%%%%%%%%%%%%%%%%%%%%%%%%%%%%%%%%%%%%%%%%%%%%%%%%%%%
\begin{num}\label{sec:lifting}  
%%%%%%%%%%%%%%%%%%%%%%%%%%%%%%%%%%%%%%%%%%%%%%%%%%%%%
Given two divisors $E_i$ and $E_j$, assume the adjacency $N_i\subseteq\overline{N}_j$ is true.
By Theorem \ref{theo:ja} there exists a wedge $\alpha:(\CC^2,O)\to X$ realising the adjacency with $\alpha_0\in\dot{N}_i$ and $\alpha_s\in\dot{N}_j$ for $s\neq 0$. Unless we state the contrary, $s$ stands for any value different from $0$ in a small enough representative $\Lambda$ of $(\CC,0)$. We consider the mapping \begin{equation}\label{eq:wedge_ext}(\alpha,pr_2):(\CC^2,O)\to X\times\CC\end{equation}
where $pr_2$ is the projection onto the second factor.

We take a Milnor radius $\epsilon$ for the image of $\alpha_0$. We denote by $B_\epsilon$ the ball in $\CC^3$ of radius $\epsilon$ and by $\SSS_\epsilon$ its boundary.  We choose a disc $\Lambda$ in $\CC$ such that $\SSS_\epsilon$ is transverse to the image of $\alpha_s$ for all $s\in \Lambda$  and such that $\alpha_s\inv(B_\epsilon)$ is a topological disc (a compact connected and simply connected neighbourhood of $0$ in $\CC$) for all $s\in \Lambda$. We denote  \begin{equation}\label{eq:calU}\calU:=(\alpha,pr_2)\inv(B_{\epsilon}\times\Lambda).\end{equation}   
We say $\underline{\alpha}:=\alpha|_{\calU}$ is \emph{a good representative of the wedge} $\alpha$. Observe that the germ  $(\underline{\alpha}\inv(O),O)$ may have other components different from $\{0\}\times\Lambda$; we call them \emph{returns}. 

Consider the situation in the following diagram:
\begin{equation}\label{dia:lifting_wedges}\xymatrix{
&&\CC^2\times \Lambda\ar_{(p,pr_2)}[dd]\\
&&\\
\calU\subseteq\CC\times\Lambda\ar^{(\underline{\alpha},pr_2)}[rr]&&X\times\Lambda
}\end{equation}
Consider the space $W^{\underline{\alpha}}:=(p,pr_2)\inv((\underline{\alpha},pr_2)(\calU))$. Then, there exists
a neighbourhood $\calV$ of $(O,0)$ in $ \CC^2\times\Lambda$, such that $W^{\underline{\alpha}}\cap \calV$ is an analytic surface. By the Principal Ideal Theorem there exists an equation for the germ $(W^{\underline{\alpha}},(O,0))$ in $\calV\subseteq \CC^2\times\Lambda$ which we denote by $G^{\underline{\alpha}}(x,y,s)=0$. 

Observe that $W^{\underline{\alpha}}$ is a family of plane curves. We denote the plane curve $W^{\underline{\alpha}}\cap(\CC\times\{s\})$ by $W^{\underline{\alpha}_s}$. 

We will need the following observation:
\begin{lem}\label{lem:reduced}
If $G^{\underline{\alpha}}$ is a reduced equation for $(W^{\underline{\alpha}},O)$, then $G^{\underline{\alpha}}_s:=G^{\underline{\alpha}}|_{\CC^2\times\{s\}}$ is reduced for all $s$ in a small enough neighbourhood $\Lambda'$ of $0$ in $\Lambda$.
\end{lem}
It is an easy consequence of the fact that the family of curves $(p,pr_2)(W^{\underline{\alpha}})\subseteq X\times\Lambda$ is parametrized in family by $(\underline{\alpha},pr_2)$, that $\underline{\alpha}_0:(\CC,0)\to X$ is one-to-one and  that $(p,pr_2):\CC^2\to X$ is a covering map outside the origin.
\end{num}

%%%%%%%%%%%%%%%%%%%%%%%%%%%%%%%%%%%%%%%%%%%%%%%%%%%%%
%\subsection{}
\begin{num}\label{sec:mult} 
%%%%%%%%%%%%%%%%%%%%%%%%%%%%%%%%%%%%%%%%%%%%%%%%%%%%%
Let $C$ be an auxiliary curve  in $\CC^2$. By abuse of notation, we denote by $C$ also its defining equation.
We study its intersection multiplicity with the curves $W^{\underline{\alpha}_s}$ in $\CC^2$. By the upper semi-continuity of the intersection multiplicity we have
\begin{equation}\label{eq:mult_wedge}%wedge1
I_O(W^{\underline{\alpha}_0},C)\geq I_O(W^{\underline{\alpha}_s},C).
\end{equation}

Let us compute first $I_O(W^{\underline{\alpha}_0},C)$. By definition of the representative $\underline{\alpha}$, the arc representative $\underline{\alpha}_0$ has only one point in $\underline{\alpha}_0\inv(0)$. Besides, since $\alpha_0$ is in $\dot{N}_i$ and because of the commutativity of the diagram (\ref{dia:indet}), we know that the strict transform of $W^{\underline{\alpha}_0}$ in $Y$ is the union of transverse curves to the dicritical divisors over $E_u$. In Table \ref{tabla1}, we have generic parametrizations of a branch of the push down of such curves in $\CC^2$. Then, it is easy to compute  $I_O(W^{\underline{\alpha}_0},C)$ for $C$ a generic line -transverse to all of the components of $(W^{\underline{\alpha}_0},O)$,- or a tangent line -corresponding with one point in $F_1$ associated to one of the dicritical divisors over $E_u$-. 

For example (see Table \ref{tabla1}), if $\alpha_0$ is in $\dot{N}_4$, we get that $W^{\underline{\alpha}_0}$ has 40 components with 20 different tangent lines corresponding to the points of $F$ (observe that each point in $E_4$ has two preimages in each dicritical divisor). Take one of these tangent lines and denote it by $L_4$. It is tangent to 2 branches and transverse to 38. Then, we have:
$$I_O(W^{\underline{\alpha}_0},L_4)=2\cdot I_O((t,\lambda t^2),y=0)+38\cdot\mult((t,\lambda t^2))=42.$$
For a generic line $C$ we have 
$$I_O(W^{\underline{\alpha}_0},C)=40\cdot\mult((t,\lambda t^2))=40.$$
We summarize these computations for any divisor $E_k$ in Table \ref{tabla2} using the following notation: 
\begin{notation} We denote by $\gamma_k$ an arc in $\dot{N}_k$ and by $\underline{\gamma}_k$ a representative of it defined in a disc $\DD$, preimage of a ball $B_{\epsilon_0}$ in $\CC^3$ with $\epsilon_0$ a Milnor radius for the image of $\gamma_k$. The curve $W^{\underline{\gamma}_k}$ is the preimage $p\inv(\underline{\gamma}_k(\DD))$ in $\CC^2$. We denote by $L_k$ one of its tangent lines. For $k\neq 1$, we also ask $\gamma_k$ being such that the strict transform of $W^{\underline{\gamma}_k}$ in $Y$ does not meet the lifting of any of its tangent lines at the origin. A parametrization of a branch of $W^{\underline{\gamma}_k}$ is given in Table \ref{tabla1}.
\end{notation}

\begin{table}[!h]%
\centering
\caption{Intersection multiplicities in $\CC^2$.} 
\label{tabla2}
\begin{tabular}{|l|l|l|l|l|l|l|l|l|}
\hline
Divisor&$E_1$&$E_2$&$E_3$&$E_4$&$E_5$&$E_6$&$E_7$&$E_8$\rule[-2.25mm]{0mm}{7mm}\\
\hline 
Number of branches of $W^{\gamma_k}$&60&30&20&40&12&24&12&24\rule[0mm]{0mm}{5mm}\\
$\mult(W^{\underline{\gamma}_k},O)$&-&30&20&40&12&24&36&48\\
$I_O(W^{\underline{\gamma}_k},L_k)$&60&32&24&42&20&30&40&50\rule[-2.25mm]{0mm}{0mm}\\
\hline
\end{tabular}

\end{table}

In particular, since $\alpha$ realises the adjacency $N_i\subseteq \overline{N}_j$, equation (\ref{eq:mult_wedge}) implies 
\begin{equation}\label{eq:mult_simple}
I_O(W^{\underline{\gamma}_i},C)\geq I_O(W^{\underline{\gamma}_j},C).
\end{equation}
By contradiction to this inequality for $C$ a transverse line to all the components, that is, comparing multiplicities of $W^{\underline{\gamma}_i}$ and $W^{\underline{\gamma}_j}$, we get the following:
\begin{equation}\label{eq:rule1}
\begin{array}{ccccccc}
  N_5\nsubseteq \overline{N}_3,& N_5\nsubseteq \overline{N}_6,&N_5\nsubseteq \overline{N}_2,& N_5\nsubseteq \overline{N}_7,&N_5\nsubseteq \overline{N}_4,&N_{5}\nsubseteq \overline{N}_{8},&N_{5}\nsubseteq \overline{N}_{1},\\
 &N_3\nsubseteq \overline{N}_6,& N_3\nsubseteq \overline{N}_2,&N_3\nsubseteq \overline{N}_7,&N_3\nsubseteq \overline{N}_4,&N_3\nsubseteq \overline{N}_8,&N_3\nsubseteq \overline{N}_1,\\
&&N_6\nsubseteq \overline{N}_2,
&N_6\nsubseteq \overline{N}_7,&N_6\nsubseteq \overline{N}_4,&N_6\nsubseteq \overline{N}_8,&N_6\nsubseteq \overline{N}_1,\\
&& &N_2\nsubseteq \overline{N}_7,&N_2\nsubseteq \overline{N}_4,&N_2\nsubseteq \overline{N}_8,& N_2\nsubseteq \overline{N}_1,\\
&&&&N_7\nsubseteq \overline{N}_4,&N_7\nsubseteq \overline{N}_8,&N_7\nsubseteq \overline{N}_1,\\
&&&&&N_4\nsubseteq \overline{N}_8,&N_4\nsubseteq \overline{N}_1,\\
&&&&&&N_8\nsubseteq \overline{N}_1.
\end{array}
\end{equation}

\end{num}
%%%%%%%%%%%%%%%%%%%%%%%%%%%%%%%%%%%%%%%%%%%%%%%%%%%%%
%\section{Proof of Nash Conjecture for $E_8$. Part II}
\subsection*{Part II}
%%%%%%%%%%%%%%%%%%%%%%%%%%%%%%%%%%%%%%%%%%%%%%%%%%%%%
%\subsection{}
\begin{num}\label{sec:riemann} 
%%%%%%%%%%%%%%%%%%%%%%%%%%%%%%%%%%%%%%%%%%%%%%%%%%%%%
Let $Y$ be the resolution of indeterminacy as in \ref{sec:indet}. For any $s\in \Lambda$ we take the arc representative $\underline{\alpha}_s:\DD\to X$ and consider the fiber product $\DD\times_{\tilde{X}} Y$ of  $\underline{\tilde{\alpha}}_s$ by $\bar{p}$. We denote it by $S^{\underline{\alpha}}_s$ as in the following diagram:
\begin{equation}\label{dia:S_gamma}\xymatrix{
S^{\underline{\alpha}_s}:=\DD\times_ {\tilde{X}} Y\ar_{q_1}[dd]\ar^{ \quad\quad\quad\quad\quad q_2}[rr]&&Y\ar_{\bar{p}}[d]\ar_{\bar{\pi}}[ddrr]&&\\
&&\tilde{X}\ar_{\pi}[d]&&\\
\DD\ar^{\underline{\alpha}_s}[rr]\ar[urr]^{\underline{\tilde{\alpha}}_s}&&X&&\CC^2\ar_{p}[ll]}
\end{equation}
It is immediate that $S^{\underline{\alpha}_s}$ is a Riemann surface and that $S^{\underline{\alpha}_s}\stackrel{\bar{\pi}\comp q_2}{\longrightarrow}W^{\underline{\alpha}_s}$ is the normalization of $W^{\underline{\alpha}_s}$. 
Since $Y$ is obtained by blowing ups in points starting in $\CC^2$, it is clear that the action of $\G$ lifts to $Y$. Then,  it also lifts to $S^{\underline{\alpha}_s}$:
\begin{eqnarray}&\hbox{ for $g\in \G$ and $(t,y)\in S^{\underline{\gamma}}\cont\DD\times Y$,}&\nonumber\\&\hbox{we define  $g(t,y):=(t,gy)$.}& \end{eqnarray}
Clearly the mapping $q_1:S^{\underline{\alpha}_s}\to \DD$ is the quotient map. It is a covering map and all the possible ramification points are in $q_1\inv(\underline{\alpha}_s\inv(O))$. Because of action of $\G$ in $S^{\underline{\alpha}_s}$ we get that all its connected components are homeomorphic. 

Assume that $\underline{\alpha}_s\inv(O)$ is only the origin. Observe that the only way of covering a disc by a connected surface only ramifying over the origin is by another disc with only one ramification point. Since $\alpha_s$ is in $\dot{N}_j$, we find that the ramification points correspond to the preimage in $Y$ of the point in $E_j$ that the lifting $\underline{\tilde{\alpha}}_s$ meets. The number of these preimages is the same that the number of branches that was given for $E_j$ in Table \ref{tabla1}, say $b_j$. Then, $S^{\underline{\alpha}_s}$ consists on $b_j$ discs. 

On the other hand, by definition of the wedge representative $\underline{\alpha}$, the image of each $\underline{\alpha}_s$ is transverse to $\SSS_\epsilon$ for any $s\in\Lambda$ and hence we get that the number of boundary components of $S^{\underline{\alpha}_0}$ and $S^{\underline{\alpha}_s}$ is always the same for all $s\in\Lambda$ with $\Lambda$ small enough. This implies that, if $b_i$ and $b_j$ (the numbers of branches given in Table \ref{tabla1}) are different, then $\underline{\alpha}_s\inv(0)$ has more points than the origin for $s\neq 0$, that is, $\underline{\alpha}$ has returns as defined right after (\ref{eq:calU}). 
\end{num}
%%%%%%%%%%%%%%%%%%%%%%%%%%%%%%%%%%%%%%%%%%%%%%%%%%%%%
%\subsection{} 
\begin{num}\label{sec:igual_bound} 
%%%%%%%%%%%%%%%%%%%%%%%%%%%%%%%%%%%%%%%%%%%%%%%%%%%%%
If $b_i=b_j$ and we assume $\underline{\alpha}$ does not have returns, then $W^{\underline{\alpha}_s}$ has the same number of components than $W^{\underline{\alpha}_0}$. In particular $W^{\underline{\alpha}}$ has $b_i$ components, denoted by $W^{\underline{\alpha},l}$ for $l=1,...,b_i$, and the intersection one of these components $W^{\underline{\alpha},l}$ with $\CC\times\{s\}$ gives a component of $W^{\underline{\alpha}_s}$ that we denote by $W^{\underline{\alpha}_s,l}$. Hence we can compare intersection multiplicities in each of these families of irreducible curves. For instance, for the adjacency $N_7\subseteq \overline{N}_5$, we have $b_7=b_5=12$. If $\underline{\alpha}$ realises this adjacency without returns, then $W^{\underline{\alpha}}$ has 12 components. For each component $W^{\underline{\alpha},l}$ we can choose the common tangent line of the curves $W^{\underline{\alpha}_s,l}$ for all $s\in\Lambda\setminus\{0\}$. We denote it by $L_5$. Then,
$$I_O(W^{\underline{\alpha}_0,l},L_5)=I_O(W^{\gamma_7,l},L_5)\geq I_O(W^{\underline{\alpha}_s,l},L_5)=I_O(W^{\gamma_5,l},L_5).$$ 
Using the parametrizations of a branch in $W^{\gamma_k}$ given in Table \ref{tabla1}, we get $I_O(W^{\gamma_5,l},L_5)=9$ and $I_O(W^{\gamma_7,l},L_5)=7$. Then the inequality is false and such a wedge is not possible. 

We reach the same conclusion for the adjacency $N_8\subseteq \overline{N}_6$. These are the only adjacencies $N_i\subseteq \overline{N}_j$ that have not been ruled out yet satisfying $b_i=b_j$. Consequently, if $\alpha$ realises an adjacency in $X$, then it has a return.
\end{num}
%%%%%%%%%%%%%%%%%%%%%%%%%%%%%%%%%%%%%%%%%%%%%%%%%%%%%
%\subsection{}
\begin{num}\label{sec:mult_gen} 
%%%%%%%%%%%%%%%%%%%%%%%%%%%%%%%%%%%%%%%%%%%%%%%%%%%%%
Assume $\underline{\alpha}$ has returns. Let $\calR$ be a component of $(\underline{\alpha}\inv(O),O)$ different from $\CC\times \Lambda$. For any $s\in\Lambda\setminus\{0\}$ small enough, the arc germs of $\underline{\alpha}_s$ centered at the points of $\calR\cap (\CC\times\{s\})$ are in the same $N_k$. We say that $\underline{\alpha}$ \emph{has a return in} $N_k$. If besides, these germs live in $\dot{N}_k$, we say that the return is \emph{transverse}. Otherwise we say it is \emph{non-transverse}. In this case, the arc germs have lifting either meeting some crossing $E_{k_1}\cap E_{k_2}$ or some $E_k$ in a non-transverse way at a smooth point.

Assume $\underline{\alpha}$ has transverse returns in $N_{k_1}$, ...,$N_{k_r}$. Then, $W^{\underline{\alpha}_s}$ for $s\neq 0$ has more branches than $W^{\gamma_j}$. In fact $W^{\underline{\alpha}_s}$ can be seen as a union of the form
$$W^{\underline{\gamma}_i}\cup W^{\underline{\gamma}_{k_1}}\cup...\cup W^{\underline{\gamma}_{k_r}}$$
for some $\gamma_l$ in $\dot{N}_l$ for $l=v,k_1,...,k_r$.

Hence, we can study the intersection multiplicity with a curve $C$ in $\CC^2$ as in paragraph \ref{sec:mult} and applying (\ref{eq:mult_wedge}) we get 
\begin{equation}\label{eq:mult_gen}%(wedge3)
  I_O(W^{\underline{\gamma}_i},C)\geq I_O(W^{\underline{\gamma}_j},C)+ \sum_{p\in \underline{\alpha}_s\inv(O)\setminus\{0\}}I_O(W^{\underline{\gamma}_{k_{p}}},C)
\end{equation}

Assume that $\underline{\alpha}$ has a return $\calR\cont \underline{\alpha}\inv(O)$ that is non-transverse. In order to upper-bound the intersection multiplicity of the branches of $(W^{\underline{\alpha}_s},O)$ that come from the return $\calR$, we take a small deformation of the arc germ of $\underline{\alpha}_s$ centered at a point in $\calR\cap(\CC\times\{s\})$: we take its lifting in $\tilde{X}$ and perturb it making it transverse to the exceptional divisor $E=\cup_{k=1}^r E_k$. If the lifting is non-transverse to a divisor $E_k$ in a generic point, then the perturbation meets at least twice $E_k$ in a transverse way. Then, the intersection multiplicity of the corresponding components with a curve $C$ will be lower bounded by $2\cdot I_O(W^{\gamma_k},C)$.
If the perturbation meets an intersection point of $E_{k_1}\cap E_{k_2}$, then the perturbation will meet transversally $E_{k_1}$ and $E_{k_2}$ and we can lower bound the intersection multiplicity of the corresponding components with a curve $C$ by $I_O(W^{\gamma_{k_1}},C)+I_O(W^{\gamma_{k_2}},C)$.

For example, if $\underline{\alpha}$ realises the adjacency $N_7\subseteq N_3$ with a non-transverse return in $N_5$, then, for instance, taking a generic line $C$ we get 
\begin{equation}\label{eq:non_trans}I_O(W^{\underline{\alpha}_0},C)=\mult(W^{\gamma_7})\geq I_O(W^{\underline{\alpha}_s},C)\geq\mult(W^{\gamma_3})+\left\{\begin{array}{l} 2\cdot\mult(W^{\gamma_5})\\ \rm{or}\\
\mult(W^{\gamma_5})+\mult(W^{\gamma_6})\end{array}\right. \end{equation}
and comparing with the values of $\mult(W^{\gamma_k})$ in Table \ref{tabla2} we get that this is false. Then such a wedge does not exists. For a transverse return in $N_5$ we can not reach a contradiction for the moment.

Since the minimum of the multiplicities $\mult(W^{\gamma_k})$ is 12, if $\alpha$ realises the adjacency $N_i\cont\overline{N}_j$ with a return (this is the only possibility by \ref{sec:igual_bound}), by inequality (\ref{eq:mult_gen}) about the intersection multiplicity with $C=L_j$ we get
\begin{equation}I_O(W^{\gamma_i},L_j)=mult(W^{\gamma_i},O)\geq I_O(W^{\gamma_j},L_j)+12.\end{equation}
Looking at Table \ref{tabla2}, we see that it is impossible for several cases. We get the following:
\begin{equation}\label{eq:rule2}\begin{array}{ccccccc}
N_1\nsubseteq \overline{N}_8,&
N_2\nsubseteq \overline{N}_3,&
N_2\nsubseteq \overline{N}_5, &
N_2\nsubseteq \overline{N}_6,&
N_3\nsubseteq \overline{N}_5,&
N_4\nsubseteq \overline{N}_2,& 
N_4\nsubseteq \overline{N}_6, \\
N_4\nsubseteq \overline{N}_7,&
N_6\nsubseteq \overline{N}_3, &
N_6\nsubseteq \overline{N}_5,&
N_7\nsubseteq \overline{N}_2,&
 N_7\nsubseteq \overline{N}_6,&
N_8\nsubseteq \overline{N}_4,& 
N_8\nsubseteq \overline{N}_7.
\end{array}\end{equation}
The remaining adjacencies are the following:
\begin{equation}\label{eq:last}\begin{array}{ll}
N_7\subseteq \overline{N}_j,\,\rm{with}\, j=3,5,&
\quad\quad N_8\subseteq \overline{N}_j,\,\rm{with}\, j=2,3,5,6,\\ N_4\subseteq \overline{N}_j,\,\rm{with}\, j=3,5,&
\quad\quad N_1\subseteq \overline{N}_j,\,\rm{with}\, j=2,3,4,5,6,7.\\
\end{array}\end{equation}
For any of these adjacencies we proceed as in the following example: the case of $\underline{\alpha}$ realising the adjacency $N_1\subseteq \overline{N}_3$:
\begin{itemize}
\item If the returns in $W^{\underline{\alpha}_s}$  have multiplicity greater than 36, then we find a contradiction with inequality (\ref{eq:mult_gen}) for $C$ a transverse line to all the branches of $W^{\underline{\alpha}_s}$. Let us study the other cases:
\item if there is a non-transverse return in $N_2$, $N_3$, $N_6$ or $N_7$ we get a contradiction as in (\ref{eq:non_trans}) with $C=L_2,\, L_3,\, L_6$ or $L_7$ respectively.
\item we can not rule out a non-transverse return in $N_5$.
\item we can not rule out 2 transverse returns in $N_5$.
\item we can not rule out transverse returns in $N_2$, $N_3$, $N_5$,  $N_6$ or $N_7$.
\item if there are 3 returns in $N_5$ or one in $N_5$ and another one in $N_6$, then 
we find a contradiction with inequality (\ref{eq:mult_gen}) for $C=L_5$.
\end{itemize}
Using the method above, in all of the 14 remaining adjacencies we can classify the possibilities for wedges realising adjacencies that have still to be ruled out. In order to avoid unnecessary repetitions, we refer the reader to their enumeration in the first column of the table \ref{tabla:strata} (for the notation there look at (\ref{eq:cases}).
\end{num}

%%%%%%%%%%%%%%%%%%%%%%%%%%%%%%%%%%%%%%%%%%%%%%%%%%%%%
%\subsection{} 
\subsection*{Part III}
\begin{num}\label{sec:versal} 
%%%%%%%%%%%%%%%%%%%%%%%%%%%%%%%%%%%%%%%%%%%%%%%%%%%%%
For any of the remaining cases, the method consists in disproving the existence of wedge representatives $\underline{\alpha}$ realising an adjacency $N_i\subseteq \overline{N}_j$ with certain returns by means of the equation of the image \begin{equation}\label{eq:wedge_image}(\alpha,pr_2)(\calU)=(p,pr_2)(W^{\underline{\alpha}})\end{equation}
 in $X\times \Lambda$. We will find in \ref{sec:implicit} all the possible equations for each separate case. 
 
By Theorem \ref{theo:ja}, it is sufficient to rule out the existence of wedges with a chosen special arc. Thus, we fix an arc $\gamma_i$ in $\dot{N}_i$ and take its reduced equation $h$ in $X$. For any wedge representative $\underline{\alpha}$ with $\alpha_0=\gamma_i$, the equation of (\ref{eq:wedge_image}) its a deformation of $h$. This implies that there is an equivalent uni-parametric deformation in the $\calR$-versal deformation of $h$ that also realises the adjacency with the same returns in the same $N_k$. 
By $\calR$-versal deformation we mean under right-equivalence by composing with biholomorphism germs $\phi$ of $X$ that admit a family of diffeomorphisms $\phi_t$ with $\phi_0=Id_X$ and $\phi_{t_0}=\phi$ for some $t_0$.

In order to study the cases in (\ref{eq:last}), we choose transverse arcs in $\dot{N}_7$, $\dot{N}_8$, $\dot{N}_4$ and $\dot{N}_1$ (see \ref{sec:implicit} for an easy way of finding them) and with the aid of \newblock {\sc Singular} compute its versal deformations: 
\begin{enumerate}
\item The versal deformation
of $y^2+z^3$, which is the equation of the image of an arc in $\dot{N}^c_7$, is
$$y^2+z^3+b_{1}x^2z+b_{2}x^2+b_{3}y^3+b_{4}y^2z+b_{5}yz^2+b_{6}xz+b_{7}y^2+$$
$$b_{8}yz+b_{9}x+b_{10}z^2+b_{11}y+b_{12}z+b_{13}.$$
\item The versal deformation of $z^4+xy$, which is the equation of the image of an arc in $\dot{N}_8^c$, is
$$z^4+xy+b_{1}x^3z+b_{2}x^3+b_{3}x^2y+b_{4}x^2z+b_{5}xy^2+b_{6}y^2z^2+b_{7}x^2+b_{8}yz^3+b_{9}xz^2+b_{10}y^2z+$$
$$b_{11}xy+b_{12}yz^2+b_{13}xz+b_{14}y^2+b_{15}z^3+b_{16}yz+b_{17}x+b_{18}z^2+b_{19}y+b_{20}z+b_{21}.$$
\item The versal deformation of $y^2+xz$, which is the equation of the image of an arc in $\dot{N}^c_4$, is
$$y^2+xz+b_{1}x^2y+b_{2}x^2z+b_{3}x^2+b_{4}yz^3+b_{5}xy+b_{6}z^4+b_{7}yz^2+$$
$$b_{8}xz+b_{9}z^3+b_{10}yz+b_{11}x+b_{12}z^2+b_{13}y+b_{14}z+b_{15}.$$
\item The versal deformation of $x^2+2y^3$, which is the equation of the image of an arc in $\dot{N}^c_1$, is
$$x^2+2y^3+b_{1}x^2yz^3+b_{2}x^2yz^2+b_{3}x^2z^3+b_{4}x^2yz+b_{5}xyz^3+b_{6}x^2z^2+b_{7}x^2y+b_{8}y^2z^3+b_{9}xyz^2+$$
$$b_{10}x^2z+b_{11}yz^4+b_{12}xz^3+b_{13}y^2z^2+b_{14}xyz+b_{15}x^2+b_{16}yz^3+b_{17}xz^2+b_{18}y^2z+b_{19}xy+$$
$$b_{20}z^4+b_{21}yz^2+b_{22}xz+b_{23}y^2+b_{24}z^3+b_{25}yz+b_{26}x+b_{27}z^2+b_{28}y+b_{29}z+b_{30}.$$
\end{enumerate}
\end{num}
%%%%%%%%%%%%%%%%%%%%%%%%%%%%%%%%%%%%%%%%%%%%%%%%%%%%%
\begin{num}\label{sec:iff}
%%%%%%%%%%%%%%%%%%%%%%%%%%%%%%%%%%%%%%%%%%%%%%%%%%%%%
Now we want to know when a uniparametric deformation of curves in X contained in any of these versal deformations is in fact the image of a wedge. We will use the following classical theorem:
\begin{theo} \label {theo:delta_cte}{\cite{Tei2}} A flat family of (space) curves admits a normalization in family if and only if it is  $\delta$-constant.
\end{theo}
In our case we prove the following::
\begin{prop}\label{prop:iff}
%Let $H_0(x,y,z)=0$ be the equation of $p(W^{\underline{\gamma}})\cont X\cont \CC^3$ for $\underline{\gamma}$ a Milnor representative of an arc. T
Consider a family of curves in $X$ given by a family of equations depending holomorphically on a complex parameter $s\in\Lambda$
$$H_s(x,y,z)=\sum_{i,j,k}a_{i,j,k}(s)x^iy^jz^k=0,\ \ a_{i,j,k}\in\CC\{s\}$$
with $H_0\inv(0)$ an irreducible curve. The family $H_s(x,y,z)=0$ is the image of a wedge representative in $X\times\Lambda$ as in (\ref{eq:wedge_image}) if and only if, given  a Milnor radius $\epsilon$ for $H_0\inv(0)$, for any $s$ in a small enough $\Lambda$ we have that 
\begin{equation}\label{eq:delta_cte}\sum_{p\in H_s\inv(0)\cap B_\epsilon}\delta(H_s\inv(0),p)=\delta(H_0\inv(0),O),\end{equation}
that is, the family ${H}_s\inv(0)\cap B_{\epsilon_0}$ is a $\delta$-constant family.
\end{prop}
\begin{proof} One direction is immediate: a family of curves given by an equation $H_s=0$ is flat; if moreover, it is the image of a wedge, then in particular it has normalization in family. By Theorem \ref{theo:delta_cte} it is $\delta$-constant. 

Let us see the other direction. If the family $H_s=0$ is $\delta$-constant, since the normalization of $H_0\inv(0)$ is a disc $\DD$ (because it is irreducible), then the normalization of the family is certain $n:\DD\times\Lambda\to H\inv(0)$. Let $\calC_1$,...,$\calC_n$ the irreducible components of the germ $(H\inv(\{0\}\times \Lambda))$ at the origin. Choose one of them, say $\calC_1$ and let $(z^m,\phi(z))$ be a Puisseux parametrization of it. We take the following mapping in $\CC^2$:
\begin{equation}\label{eq:Phi}\Phi:(t,z)\mapsto (z^m,\phi(z)+t).\end{equation}   
Then, $\alpha:=n\comp \Phi:\DD\times\Lambda\to X$ has the same image than $n$ and it is a wedge: $\alpha\inv(O)$ contains $\{0\}\times\Lambda$.
\end{proof}

 \end{num}
%%%%%%%%%%%%%%%%%%%%%%%%%%%%%%%5
%\subsection{}
\begin{num}\label{sec:non-trans}
%%%%%%%%%%%%%%%%%%%%%%%%%%%%%%%%%%%%%%%%%%%%%%%%%%%%%
 Among the adjacencies in (\ref{eq:last}), there are only 3 possible cases of wedges realising adjacencies with non-transverse returns that can not be ruled out in the way of (\ref{eq:non_trans}): a wedge realising $N_1\subseteq\overline{N}_3$ with a non-transverse return in $N_5$ or a wedge realising $N_1\subseteq\overline{N}_5$ with non-transverse return either in $N_3$ or in $N_5$. We can be more precise: the lifting of the arc germs associated to these non-transverse returns, meet the corresponding divisor in at a smooth point and with intersection multiplicity at this point equal to $2$. % have lifting to $\tilde{X}$ non-transverse through a smooth point in the corresponding divisor. %These equations are in the closure of the corresponding ones with transverse returns in the following sense:
In these cases we have the following, proved in a classical way:
\begin{prop}\label{prop:non_trans}If there is a wedge $\alpha$ realising the adjacency ${N}_i\subseteq\overline{N}_j$ in any of these cases with a non-transverse return in $N_k$ , then there is a wedge realising the adjacency only with transverse returns, all in $\dot{N}_k$. 
\end{prop}
\begin{proof} Let $H_s(x,y,z)=0$ be the equation of the image of $(\underline{\alpha},pr_\Lambda)$ in $X\times\Lambda$. After a possible base change of type (\ref{eq:Phi}) we can assume that the singularities in $H_s\inv(0)$ out of the origin are parametrized, say by $p_l(s)$ with $l=1,...,n$, for $p_l:(\Lambda,0)\to (X,O)$ with $pr_2\comp p_l=id_\Lambda$.

%In the same way than in the proof of Theorem \ref{theo:delta_cte_E_8}, 
We claim that we can construct a deformation of $h_0$ with base $\Lambda\times\CC^N$ of the form 
\begin{equation}\label{eq:family}G(x,y,z,s,\underline{v})=H_{s}(x,y,z)+Q(x,y,z,s,\underline{v}),\end{equation}
with $Q$ a family of polynomials in three variables $x,y,z$ depending on $s$ and on other parameters denoted by $\underline{v}\in\CC^N$, (that is, $Q$ belongs to $\CC\{s,\underline{v}\}[x,y,z]$), such that the deformation $G$ satisfies the following equalities 
$$G(x,y,z,0,\underline{v})=H_0(x,y,z)\quad {\rm and}\quad    G(x,y,z,s,0,...,0)=H_s(x,y,z)$$
for any $\underline{v}$ and $s$ respectively, and the following additional properties:
\begin{enumerate}
%\item the germ of $G_{t,v_1,...,v_N}$ at the origin has a 
\item The deformation is versal at  $(0,0,0,s,0,...,0)$ for any $s\neq 0$.
\item  the germ of $G_{s,\underline{v}}$ at $p_l(s)$
has a singularity topologically equivalent to the one that $H_s$ has at the same point for any $l\leq n$, any $s\neq 0$ and any $\underline{v}$.
\end{enumerate}
Assuming we can construct such a deformation, we end the proof as follows. 

Define $\calZ\subset \Lambda\times\CC^N$ to be the stratum of the base of the deformation $G$ consisting of points $b$ such that $G_b=0$ has the topological type of the curve associated to the transverse case. %This a semi-analitycal subset.   
The origin is at the closure of $\calZ$ since   $(\Lambda\setminus\{0\})\times\{O\}$ is at the closure of $\calZ$. This is because we can take the following deformation of $\underline{\alpha}_s\inv(0)$ given by a deformation of the parametrizations of its branches: leaving all parametrizations fixed except the one corresponding to the non-transverse return that we perturb in order to get curves with lifting transverse to $E_k$ (we perturb the non-transverse lifting and blow-down). Observe that this is a $\delta$-constant deformation because it is a deformation of parametrizations. %It is a flat deformation because it is a deformation of parametrizations and hence %\textbf{E}_8$ is a UFD and it is a $\delta$-constant deformation because of the definition of $\gamma$ and (ii). 

By Curve Selection Lemma there exists an arc $\gamma:(\CC,0)\to \overline{\calZ}$ such that $\gamma(s)\in\calZ$ for all $s\neq 0$ and $\gamma(0)=O$. We take the deformation $G_{\gamma(t)}$. 
%Since $G$ is versal at $(0,0,0,s,0,...,0)$ by the property (i), this deformation is equivalent to a deformation pull back of $G$: 
%there exists an arc $\gamma:(\CC,0)\to (X\times \Lambda\times \CC^N,O)$ such that $G\comp \gamma$ is a $\delta$-constant deformation of $H_0$ with generic curve in $\dot{N}_j+\dot{N}_k+\dot{N}_k+...+\dot{N}_k$.
%The deformation is $\delta$-constant as a consequence of (ii) and the definition of $\gamma$. 
It is a $\delta$-constant deformation because of the definition of $\gamma$ and (ii). 
Using Proposition \ref{prop:iff}, we get a wedge with the desired properties. 

Now we prove the claim. In order to prove the existence of $G$, we  find the appropriate $Q$. 

 Let $\mm$ denote the maximal ideal of
$\calO_{X,O}$. %Choose $K$ big enough so that, for any $s\neq 0$, the $K$-th power of $\mm$ is contained in the Jacobian ideal of $H_s$ at the origin. 
Let $R$ be the whole family of polynomials of degree at most $K$ in three variables:
$$R=\sum_{i_1+i_2+i_3\leq K}v_{i_1,i_2,i_3}x^{i_1}y^{i_2}z^{i_3}$$
with the coefficients $v_{i_1,i_2,i_3}$ varying in $\CC$; we denote the set of coefficients by $\underline{v}$. Choose $K$ big enough so that $H_{s}(x,y,z)+R(x,y,z,\underline{v})$ is versal at the origin $(0,0,0,s,0...,0)$.
%For any $s\neq 0$, the deformation of $H_{s}$ given by $H_{s}(x,y,z)+R(x,y,z,\underline{v})$ is versal at the origin $(0,0,0,s,0,...,0)$.

We look for $Q$ in the family of polynomials of degree at most a fixed $K$.

The general such polynomial may be written as
$$Q=\sum_{i_1+i_2+i_3\leq K}C_{i_1,i_2,i_3}x^{i_1}y^{i_2}z^{i_3}$$ and we denote by $\underline{C}_I$ the set of coefficients.

For any $R$, we look for $Q$ satisfying (in order that $G$, defined by (\ref{eq:family}), paralelly satisfies the previous requirements (i)-(ii)):
\begin{enumerate}
\item the $K$-Taylor expansion at the origin $(0,0,0,s,0,...,0)$ coincides with $R$.

\item the $K_l$-th Taylor expansion at $p_l(t)$ vanishes for any $l\leq n$ where $K_l$ is  the $R$-equivalence determinacy degrees of the
singularities that $H_{s}$ has at $p_l(s)$;
\end{enumerate}

These conditions give linear conditions in the unknowns $\underline{C}_I$. The matrix of coefficients of the linear system has entries in the ring $\CC\{s\}$ and the column of independent terms is formed by zeroes and the ${v}_I$'s.
It can be checked that, if $K$ is big enough, then the matrix will have maximal rank for $s\neq 0$ big enough. %Hence, there is a maximal minor of the coefficient matrix whose determinant, that depends only on $t$, does not vanish identically for $t$ small enough, and we can assume that $0$ is the only zero of the determinant. 
Solving the system by Cramer's rule with respect to the
given maximal minor we find 
$$Q(t,\underline{v})=\sum_{i_1+i_2+i_3\leq K}C_{i_1,i_2,i_3}(s,\underline{v})x^{i_1}y^{i_2}z^{i_3}$$
with coefficients $C_{i_1,i_2,i_3}$ belonging to the ring of Laurent convergent power series $\CC\{s\}[s^{-1}]$ and depending linearly on $v_I$. 

Choose $A$ a sufficiently big integer so that $s^AC_{i_1,i_2,i_3}(s,\underline{v}_I)$ are holomorphic and vanish at $0$ for any $i_1+i_2+i_3\leq N$. Then, redefining  $Q:=s^A\cdot Q$ and using formula (\ref{eq:family}) we have the required deformation $G$.

\end{proof}
\end{num}
%%%%%%%%%%%%%%%%%%%%%%%%%%%%%%%%%%%%%%%%%%%%%%%%%%%%%
%\subsection{}
\begin{num}\label{sec:implicit} 
%%%%%%%%%%%%%%%%%%%%%%%%%%%%%%%%%%%%%%%%%%%%%%%%%%%%%
We find the reduced equations for the cases of wedges with transverse returns that remain to be analysed. It will become clear that the equations in $X\times\Lambda$ do not differentiate between the images for possible wedges realising the adjacency ${N}_i\subseteq\overline{N}_j$  with a transverse return in $N_k$ or a wedge realising  ${N}_i\subseteq\overline{N}_k$ with a transverse return in $N_j$. Thus, we will denote such a case by: 
\begin{equation}\label{eq:cases}(N_i;\dot{N}_j+\dot{N}_k).\end{equation}

We find first the possible reduced equations for  $W^{\underline{\alpha}}$ in $\CC^2\times\Lambda$. Since it is a $\G$-invariant set, its reduced equation has character $1$ and then it can be expressed as a power series in  $E$, $F$ and $V$ with coefficients depending holomorphically on $s$. Then, performing the substitution given by (\ref{eq:p}), that is 
\begin{eqnarray}\label{eq:EFV}x=E,\,\,\,\,& y=F,&\,\, \,\,z=V,\end{eqnarray} 
 we get the reduced equations of (\ref{eq:wedge_image}), a power series in $x$, $y$ and $z$ with coefficients depending  holomorphically on $s$.
 
For a wedge representative $\underline{\alpha}$ with certain return, the curve germ $(W^{\underline{\alpha}_s},O)$, for $s\neq 0$, is a $\G$-invariant curve. We will denote by $\Gamma$ a $\G$-invariant curve in $\CC^2$. We will use the following notation: 
\begin{itemize}
\item we say $\Gamma\in \dot{N}_{k}$ if the strict transform of $p(\Gamma)\cont X$ by the resolution map $\pi$ meets transversally $E_{k}$ at a non-singular point. Equivalently we can ask the strict transform of $\Gamma$ in $Y$ to meet transversally the dicritical divisors over $E_{k}$ at non-singular points.
\item we say $\Gamma\in N_k^{nt[m]}$ if  the strict transform $\widetilde{p(\Gamma)}$ of $p(\Gamma)\cont X$ by the resolution map $\pi$ meets non-transversally $E_{k}$ with multiplicity $m$, that is $I_O(\widetilde{p(\Gamma)},E_k)=m$.
\item we say $\Gamma\in N_{k_1}\cap N_{k_2}$ if $\widetilde{p(\Gamma)}$ meets $E_{k_1}\cap E_{k_2}$.
\item we say that $\Gamma$ belongs to a sum of these cases, if $p(\Gamma)$ has several components, each one in one of the cases. 
\end{itemize}
If $\Gamma$ is in $\dot{N}_{k_1}+...+\dot{N}_{k_r}$, then its germ at the origin is of the form $W^{\underline{\gamma}_{k_1}}\cup ...\cup W^{\underline{\gamma}_{k_r}}$ for certain $\gamma_{k_i}\in\dot{N}_{k_i}$ for $i=1,,,.r$. In particular, in Table \ref{tabla1} we can find parametrizations of a branch of the $W^{\gamma_k}$.  

Since $\underline{\alpha}_s$ for the cases that we have to study never has the arc germ $\alpha_s$ or the returns in $N_1$, we know that the tangent cone of $W^{\underline{\alpha}_s}$ is a product of $E$, $F$ and $V$. Therefore it will be sufficient to classify $\G$-invariant curves with these tangent cones.  In fact, such an invariant curve will be the union of $\G$-invariant curves with tangent cones of the form $E^k$, $F^k$ or $V^k$ and its equation can be obtained as the product of the equations these simpler ones. 

We study separately the equations of $\G$-invariant curves $\Gamma$ with tangent cone, denoted by $TC$, equal to $E^{a}$, $F^a$ or $V^a$. We write the equations for $\Gamma$, in terms of $E$, $F$ and $V$ for the needed cases. We only write the homogeneous forms of low degree, that is the first ``monomials'' in $E$, $F$, and $V$ of the equations with respect to the weights $(30,20,12)$: any other monomial in $E$, $F$, and $V$ with smaller degree than any of those written does not appear.

\begin{enumerate}
%%%%%%%%%%%%%%%%%%%
\item[$TC=E.^{\ }$] There is only one possibility for $\Gamma$: being in $\dot{N}_2$. If $\Gamma$ is in  $N_2^{nt}$ or in $N_2\cap N_1$, then $\Gamma$ would have tangent cone with more multiplicity, as we can see bounding below the  multiplicity by means of a generic deformation of the curves which cut at least either twice $E_2$ or $E_1$ and $E_2$. The equation for a $\G$-invariant curve $\Gamma$ in the case $\dot{N}_2$, the only case with $TC=E$ is:
\begin{eqnarray}\label{eq:gen1}\Gamma\in\dot{N}_2:\,\,\,\,&a_1E+a_2FV+...&\,\,\,\,a_1\neq 0.\end{eqnarray}
%%%%%%%%%%%%%%%%%%%
\item[$TC=V.^{\ }$] Reasoning like for $TC=E$, there is only one possibility: being in $\dot{N}_5$.
\begin{eqnarray}\Gamma\in\dot{N}_5:\,\,\,\,&a_1V+a_2F+...&\,\,\,\,a_1\neq 0.\end{eqnarray}
%%%%%%%%%%%%%%%%%%%
\item[$TC=F.^{\ }$] The only possibility is being in $\dot{N}_3$.
\begin{eqnarray}\Gamma\in\dot{N}_3:\,\,\,\,&a_1F+a_2V^2+...&\,\,\,\, a_1\neq 0.\end{eqnarray}
%%%%%%%%%%%%%%%%%%%
\item[$TC=F^2$.] There are $3$ main possibilities for $\Gamma$: being in $\dot{N}_4$ or in  $\dot{N}_3+\dot{N}_3$. 
The cases $N_3\cap N_4$ and ${N}_4+{N}_1$ are ruled out by bounding below multiplicity using a generic deformation as in the case $TC=E$.
The case of $N^{nt[2]}_3$ is in the closure of $\dot{N}_3+\dot{N}_3$. 

We compute intersection multiplicity of $\Gamma$ in these cases with a generic arc $(t,\lambda t^5)$ in $\dot{N}_3$. We use the parametrizations of the different branches given in Table \ref{tabla1}. We get that for $\dot{N}_4$ it is $42$, while for $\dot{N}_3+\dot{N}_3$ is $48$.
This means that $EV$ has coefficient non-zero for the first case and $0$ for the other.
\begin{eqnarray}\Gamma\in\dot{N}_4:\,\,\,\,&a_1F^2+a_2EV+...&\,\,\,\, a_1,a_2\neq 0\\
\Gamma\in\dot{N}_3+\dot{N}_3:\,\,\,\,&a_1F^2+a_2FV^2+...&\,\,\,\, a_1\neq 0, \,\,a_i\ \rm{generic}\end{eqnarray}
%%%%%%%%%%%%%%%%%%%%%%%%%%%%%%%%
\item[$TC=V^2$.] There are $2$ main possibilities for $\Gamma$: in $\dot{N}_6$ or in $\dot{N}_5+\dot{N}_5$. Other cases as $N_5\cap N_6$ or ${N}_6\cap{N}_7$ are ruled out by bounding below multiplicity using a generic deformation. The case $N_5^{nt[2]}$ is in the closure of the case $\dot{N}_5+\dot{N}_5$.

We compute intersection multiplicity of $\Gamma$ with generic curve $(t,\lambda t^9)$ in $\dot{N}_5$ looking at the parametrizations of the different branches in Table \ref{tabla1}. We get that for $\dot{N}_6$ it is $30$, while for $\dot{N}_5+\dot{N}_5$ is $40$.% for the other case (observe that the tangent line $L_3$ is in $\dot{N}_3$). 
 This implies that for the first case the coefficient of $E$ is non-zero and for the second one it is zero:
\begin{eqnarray}\Gamma\in\dot{N}_6:\,\,\,\,&a_1V^2+a_2E+...&\,\,\,\,a_1,a_2\neq 0\\
\Gamma\in\dot{N}_5+\dot{N}_5:\,\,\,\,&a_1V^2+a_2FV+...&\,\,\,\,a_1\neq 0,\ \ a_i\  \rm{generic}\end{eqnarray}
%%%%%%%%%%%%%%%%%%%%%%%%%%%%%

\item[$TC=V^3$.] There are $3$ main possibilities for $\Gamma$: in $\dot{N}_7$, $\dot{N}_5+\dot{N}_6$ and $\dot{N}_5+\dot{N}_5+\dot{N}_5$. The cases of non-transverse curves in $N_5$ are in the closure of $\dot{N}_5+\dot{N}_5+\dot{N}_5$ and the case of $N_5\cap N_6$ is in the closure of $\dot{N}_5+\dot{N}_6$. 

We compute the intersection multiplicity with a generic curve $(t,\lambda t^9)$ in $\dot{N}_5$:
\begin{equation}\begin{array}{rcll}
I_o(\Gamma,(t,\lambda t^9))&=&40&\rm{with }\,\, \Gamma\in\dot{N}_7\\

I_o(\Gamma,(t,\lambda t^9))&=&50&\rm{with }\,\, \Gamma\in\dot{N}_5+\dot{N}_6\\

I_o(\Gamma,(t,\lambda t^9))&=&50&\rm{with }\,\, \Gamma\in N_5\cap N_6\,\,\rm{with}\,\,tangent\,\, cone\ V^3\\

I_o(\Gamma,(t,\lambda t^9))&=&60&\rm{with }\,\, \Gamma\in\dot{N}_5+\dot{N}_5+\dot{N}_5\\

I_o(\Gamma,(t,\lambda t^9))&=&60& \rm{with }\,\, \Gamma\in\dot{N}_5+ N_5^{nt[2]}\\

I_o(\Gamma,(t,\lambda t^9))&=&60&\rm{with }\,\, \Gamma\in N_5^{nt[3]}.
\end{array}\end{equation}
%and $W^{\underline{\gamma}_{i_1}}\cup...\cup W^{\underline{\gamma}_{i_k}}$ for $\underline{\gamma}$ in the corresponding cases (and always with tangent cone $V^3$).  
For the transverse cases, we compute it using the parametrizations of the irreducible components of $\Gamma$ in Table \ref{tabla1}. For the non-transverse cases we see that this intersection multiplicity is the same that in the transverse cases. For example, assume that   $\Gamma\in N_5\cap N_6$ (with tangent cone $V^3$). We consider a parametrization of the curve, which has only one component. We perturb its lifting to $\tilde{X}$ getting a parametrized curve meeting $E_5$ and $E_6$ transversally. 
Since the intersection multiplicity with the divisors $E_5$ and $E_6$ is conservative, either the intersection multiplicity is the same or the total number of intersection points of the perturbation with $E_5\cup E_6$ is at least 3. The second possibility is impossible because the tangent cone of the deformation would have more multiplicity than $V^3$. 

In order to get the final equations, we look at the intersection multiplicities of the test curve $(t,\lambda t^9)$ with the forms $E$, $F$ and $V$:
\begin{table}[!h]%
\label{fig:aux_V^3}
\centering
\begin{tabular}{|c||c|c|c|}
\hline
Inv. Form $\Sigma$&$E$&$F$&$V$\rule[-2.25mm]{0mm}{7mm}\\
\hline
\rule[-2.25mm]{0mm}{5mm}
$I_o(\Sigma,(t,\lambda t^9))$&$30$&$20$&$20$\\
\hline
\end{tabular}
\end{table}

Finally we get:
\begin{eqnarray}\Gamma\in\dot{N}_7:\,\,\,\,&a_1V^3+a_2F^2+...&\,\,\,\,a_1,a_2\neq 0\\
\Gamma\in\dot{N}_5+\dot{N}_6:\,\,\,\,&a_1V^3+a_2EV+...
&\,\,\,\,a_1\neq 0,\ a_i\ \rm{generic}\\
%N_5\cap N_6:&V^3+EF+F^2V+...&a_1\neq 0\\
\Gamma\in\dot{N}_5+\dot{N}_5+\dot{N}_5:\,\,\,\,&a_1V^3+a_2FV^2+a_3V^4&\nonumber\\
&+a_4F^2V+...&\,\,\,\,a_1\neq 0,\ a_i\ \rm{generic}
\end{eqnarray}
%More precisely, in the second case, generic means $a_2$ or $a_5\neq 0$.

%%%%%%%%%%%%%%%%%%%%%%%%%%%
\item[$TC=V^4$.]
There are 5 main possibilities for $\Gamma$: in $\dot{N}_8$, $\dot{N}_7+\dot{N}_5$, $\dot{N}_6+\dot{N}_6$,  $\dot{N}_5+\dot{N}_5+\dot{N}_6$,  %$\dot{N}_5+{N}_5\cap N_6$ 
and $\dot{N}_5+\dot{N}_5+\dot{N}_5+\dot{N}_5$. 
Other possible cases include non-transverse components to $E_5$ which are in the closure of $4\dot{N}_5$ or $\dot{N}_5+\dot{N}_5+\dot{N}_6$, non-transverse curves to $E_6$, in the closure of $\dot{N}_6+\dot{N}_6$ and the case of $N_5+N_5\cap N_6$, which is in the closure of $\dot{N}_5+\dot{N}_5+\dot{N}_6$.

Intersecting with a generic arc in $\dot{N}_5$, say $(t,\lambda t^9)$, and with with another arc $(t,\eta t^4)$ (with $\eta$ also generic) we  distinguish all the cases.  The results are summarized in the following table: 
\begin{table}[!h]%
\label{fig:Possibilities for tangent cone $V^4$.}
\centering
\begin{tabular}{|r||c|c|c|c|}
\hline

&$\dot{N}_8$&$\dot{N}_5+\dot{N}_7$&$2\dot{N}_6$&$N_6^{non-trans}$\rule[-2.25mm]{0mm}{7mm}\\
\hline
$I_O((t,\lambda t^9),\Gamma)$&$50$&$60$&$60$&$60$\rule[0mm]{0mm}{5mm}\\
$I_O((t,\lambda t^4),\Gamma)$&&$55$&$60$&$60$\rule[-2.25mm]{0mm}{0mm}\\
\hline
%\end{tabular}
%\end{table}
%\begin{table}[!h]%
%\label{fig:$V^4$.}
%\caption{Curves in $\CC^2$ with tangent cone $V^4$.}
%\begin{tabular}{|r||c|c|c|c|}
\hline
%&&&&\\
\hline
&$2\dot{N}_5+\dot{N}_6$&$\dot{N}_6+N_5^{non-trans}$&$4\dot{N}_5$&$N_5^{non-trans}+..$\rule[-2.25mm]{0mm}{7mm}\\
\hline
$I_O((t,\lambda t^9),\Gamma)$&$70$&$70$&$80$&$80$\rule[-2.25mm]{0mm}{7mm}\\
%$I_O((t,\lambda t^4),W^{\gamma})$&&&&$80$\\
\hline
%\rule[0mm]{0mm}{4mm}
\end{tabular}
\end{table}

In order to get the final equations, we look at the intersection multiplicities of the test curves $(t,\lambda t^9)$ and $(t,\lambda t^4)$ with the forms $E$, $F$ and $V$: %which can be read in the last columns of Table \ref{fig:}.
\begin{table}[!h]%
\label{fig:Possibilities for tangent cone $V^4$.}
\centering
\begin{tabular}{|r||c|c|c|}
\hline
\rule[-2.25mm]{0mm}{7mm}
Inv. Form $\Sigma$&$E$&$F$&$V$\\
\hline
\rule[0mm]{0mm}{5mm}
$I_O((t,\lambda t^9),\Sigma)$&$30$&$20$&$20$\\
$I_O((t,\lambda t^4),\Sigma)$&$30$&$20$&$15$\\
\hline
\end{tabular}
\end{table}

Finally we get:
\begin{eqnarray}\Gamma\in\dot{N}_8:\,\,\,\,&a_1V^4+a_2EF+...&\,\,\,\,a_1,a_2\neq 0,\\
\label{eq:N_5+N_7}
\Gamma\in\dot{N}_5+\dot{N}_7:\,\,\,\,&a_1V^4+a_2F^2V+...&\,\,\,\,a_1\neq 0, \,\,a_i\ \rm{generic} \\
\Gamma\in2*\dot{N}_6:\,\,\,\,&a_1V^4+a_2EV^2+...&\,\,\,\,a_1\neq 0,\,\, a_i\ \rm{generic} \\
%\dot{N}_5+N_5\cap N_6:&V^4+...&\\
\Gamma\in2*\dot{N}_5+\dot{N}_6:\,\,\,\,&a_1V^4+a_2EV^2+a_3FV^3&\nonumber\\
&+a_4V^5+a_6EFV+...&\,\,\,\,a_1\neq 0, \,\,a_i\ \rm{generic} \\
\Gamma\in4*\dot{N}_5:\,\,\,\,&a_1V^4+a_2FV^3+a_3V^5+a_4F^2V^2&\nonumber\\
&+a_5EV^3+a_6FV^4+a_7V^6+...&\,\,\,\,a_1\neq 0, \,\,a_i\ \rm{generic} 
\end{eqnarray}

%In the second case generc means that the coefficeint of $E^2$ is non-zero. In the third case, generic means $a_4-a_5\neq 0$. 
\end{enumerate}

Multiplying the equations obtained for the groups of components with the same reduced tangent cone $E$, $F$ or $V$, we get the equations for the generic arcs in other cases, for instance:
\begin{eqnarray}
%\dot{N}_2+\dot{N}_3:&a_1EF+...&a_1\neq 0\\
%\dot{N}_2+\dot{N}_5:&a_1EV+...&a_1\neq 0\\
%\dot{N}_3+\dot{N}_5:&a_1VF+...&a_1\neq 0\\
%\dot{N}_4+\dot{N}_5:&a_1VF^2+a_2EV^2+...&a_1,a_2\neq 0\\
%\dot{N}_5+2*\dot{N}_3:&a_1VF^2+a_2FV^3+...&a_1\neq 0,\ a_i\  \rm{generic}\\
\Gamma\in\dot{N}_3+\dot{N}_6:\,\,\,\,& a_1FV^2+a_2V^4+a_3FE+...&\,\,\,\,a_1,a_3\neq 0.%\\
%\dot{N}_3+2*\dot{N}_5:&a_1FV^2+a_2V^4+a_3F^2V+...&a_1\neq 0,\ a_i\ \rm{generic}\\
%\dot{N}_2+\dot{N}_6:&a_1EV^2+a_2FV^3+a_3E^2+...&a_1,a_3\neq 0\\
%\label{eq:gen2}\dot{N}_3+\dot{N}_7:&a_1FV^3+a_2F^3+...&a_1,a_2\neq 0
\end{eqnarray}
Now, in order to find the equations for $W^{\underline{\alpha}}$ with $\alpha$ a wedge realising an adjacency with certain returns, we just have to sum the equation of the special arc plus the possible equations for $p(\Gamma)=p(W^{\underline{\alpha}_s})$ with $\Gamma$ in the corresponding case. The equation of $p(\Gamma)$ is obtained from the one of $\Gamma$ after the change of variables (\ref{eq:EFV}). For instance, for a wedge $\alpha$ realising $N_{1}\subseteq N_3$ with a transverse return in $N_6$, taking the equation $x^2+2y^3$ for the special arc as in \ref{sec:versal}, we see that the equation of $W^{\underline{\alpha}}$ is of the form 
$$H(x,y,z,s)=x^2+2y^3+a_1yz^2+a_2z^4+a_3xy+....\,\,\rm{with}\,\,a_1,a_3\neq 0,$$
written in the weighted order $(30,20,12)$.
We can apply this information to find the strata of the base of the versal deformation of $x^2+2y^3$ of curves with lifting transverse to $E_3$ and $E_6$ as we will do in Table \ref{tabla:strata}.
\end{num}
%%%%%%%%%%%%%%%%%%%%%%%%%%%%%%%%%%%%%%%%%%%%%%%%%%%%%
%\subsection{}
\begin{num}\label{sec:codim} 
%%%%%%%%%%%%%%%%%%%%%%%%%%%%%%%%%%%%%%%%%%%%%%%%%%%%%
We study  condition (\ref{eq:delta_cte}) in the families we have found. We recall the following notion:

\begin{defn} Let $h_0$ be the equation of a curve in $X$. Let $H:X\times\calB\to \CC$ be its versal deformation as above and $\calT$ a topological type (a $\mu$-class) adjacent to $h_0$. Let $S_\calT$ be the stratum of $\calB$ consisting of points $b$ such that $H_b\inv(0)$ has topological type $\calT$ at the origin.
Define $\Delta=\delta(h_0)-\delta(\calT)$. In this context, the  $(\calT,\Delta)$-\emph{constant stratum} is defined to be the subset of $S_\calT$ consisting of points $b$ such that the sum of $\delta$-invariants of points in $H_b\inv(0)\setminus\ori$ is equal to $\Delta$. We denote it by $S_{\calT,\Delta}$.

$$S_{\calT,\Delta}:=\{b\in S_{\calT}:
\sum_{p\in H_b\inv(0)\cap B_{\epsilon_0}\setminus\{O\}}\delta(H_b\inv(0),p)=\Delta\}$$%(N_i;N_j+N_{i_1}+...+N_{i_k})\}.$$
%When we call \emph{the codimension of} $S_{\calT,\Delta}$, we mean  the codimension of $S_{\calT,\Delta}$ in $S_{\calT}$.
\end{defn}

Given a possible adjacency $N_i\subseteq \overline{N}_j$, we denote each case $(\dot{N}_i;\dot{N}_j+\dot{N}_{k_1}+\dot{N}_{k_2})$ by $I=(k_1,k_2)$. For each case $I$, let $\calT_I$ be the topological type of the germ $(\underline{\alpha}_s(\DD),O)$. We compute $\Delta_{i,I}=\delta(h_0)-\delta(\calT_I)$. Then we have the following result that will lead us to finish the proof:
\begin{prop}\label{prop:codim} Let $N_i\subseteq \overline{N}_j$ be an adjacency that remains to be ruled out. Let $I_1$,...,$I_m$ be the cases $(\dot{N}_i;\dot{N}_j+\dot{N}_I)$ that still remain to be ruled out (where $\dot{N}_I$ stands for $\dot{N}_{k_1}+\dot{N}_{k_2}$ for $I=\{k_1,k_2\}$). Let %$\calT_I$ the $\mu$-constant class as in \ref{eq:S_calT for the case $I$ and 
$\Delta_{i,I}=\delta(0)-\delta(\calT_I)$. In this context, if \begin{equation}\label{eq:ineq_codim}\codim(S_{\calT_I},S_{\calT_I;\Delta_{i,I}})\geq \codim(S_{\calT_I},\calA_i)\end{equation} for every case $I$, then $N_i\nsubseteq\overline{N}_j$.
\end{prop}

\begin{proof} By Proposition \ref{prop:non_trans}, in order to rule out an adjacency, it is enough to rule out the transverse cases  $(N_i;\dot{N}_j+\dot{N}_I)$.

If $\codim_{S_{\calT_I}}(S_{\calT_I;\Delta_{i,I}})\geq \codim_{S_{\calT_I}}(\calA)$
for any transverse case $I$, then it is clear that there exists  
$b_0\in\calA\setminus\bigcup_I\overline{S_{\calT_I;\Delta_{i,I}}}$. We can assume that $H$ is still versal at $b_0$ by openness of versality.  Then, since $b_0\notin\bigcup_I\overline{S_{\calT_I,\Delta_I}}$  there is no wedge realising the adjacency with $H_{b_0}$ as image of the special arc.
By Theorem \ref{theo:ja}, if there is no wedge realising the adjacency with $H_{b_0}$ as image of the special arc, then %there is no wedge realising the adjacency and
the adjacency is impossible.
\end{proof}
\end{num}
%%%%%%%%%%%%%%%%%%%%%%%%%%%%%%%%%%%%%%%%%%%%%%%%%%%%%
%\subsection{}
\begin{num}
%%%%%%%%%%%%%%%%%%%%%%%%%%%%%%%%%%%%%%%%%%%%%%%%%%%%%
Given any of the remaining cases $(N_i;\dot{N}_j+\dot{N}_{k_1}+\dot{N}_{k_2})$, we consider the versal deformation of the equation of the image of an arc in $\dot{N}_i$ given in (\ref{sec:versal}). Knowing the form of the equations of curves associated to the cases $N_i$ and  $\dot{N}_j+\dot{N}_{k_1}+\dot{N}_{k_2}$ it is easy to get the equations of the strata $\calA_i$ and $S_{\calT_I}$. In particular, the strata of the special arcs in the corresponding versal deformations are the following:
\begin{itemize}
\item $\calA_7=\{b_{8}=...=b_{13}=0\}$ in the versal deformation (i)
\item $\calA_8=\{b_{12}=...=b_{21}=0\}$ in  versal deformation (ii)
\item $\calA_4=\{b_{9}=...=b_{15}=0\}$ in the versal deformation (iii)
\item $\calA_1=\{b_{16}=...=b_{30}=0\}$ in versal deformation (iv)
\end{itemize}
We write the equations for $S_{\calT_I}$ in the base of the versal deformation and of $\calA_i$ in $S_{\calT_I}$ in Table \ref{tabla:strata}. The codimension $\codim(S_{\calT_I},\calA_i)$ is then immediate from the equations of the strata. 

We observe by looking at the equations that in all these cases
\begin{equation}\label{eq:codim_calA} \codim(S_{\calT_I},\calA_i)=\Delta_{i,I}.
\end{equation}

\begin{table}[!h]%
\centering
\caption{Strata in the deformation of the equation of a transverse arc.} 
\label{tabla:strata}
\begin{tabular}{|l|c|c|c|}
\hline
Case&$\Delta_{i,I}$&Equations of $S_{\calT_I}$&Equations of $\calA_i$\rule[-2.25mm]{0mm}{7mm}\\
\hline
$(N_7, \dot{N}_3+\dot{N}_5)$ & $1$ &$b_{9}=...=b_{13}=0$& $b_{8}=0$\rule[0mm]{0mm}{5mm}\\
%____________________________________________________ 
$(N_7, \dot{N}_5+\dot{N}_5)$ & $2$ &$b_{9}=b_{11}=...=b_{13}=0$& $b_{8}=b_{10}=0$\rule[-2.25mm]{0mm}{0mm}\\
\hline 
%____________________________________________________
 $(N_8,\dot{N}_2+\dot{N}_5)$ & $2$ &$b_{14}=...=b_{21}=0$&$b_{12}=b_{13}=0$\rule[0mm]{0mm}{5mm}\\
%________________________________________________
 $(N_8,\dot{N}_3+\dot{N}_3)$ & $2$&$b_{11}=b_{13}=...=b_{21}=0$& $b_{12}=b_{14}=0$\\
%____________________________________________________
 $(N_8,\dot{N}_3+\dot{N}_5)$ & $5$ &$b_{17}=...=b_{21}=0$& $b_{12}=...=b_{16}=0$\\
%____________________________________________________
 $(N_8,\dot{N}_3+\dot{N}_6)$ & $1$ &$b_{13}=...=b_{21}=0$& $b_{12}=0$\\
%____________________________________________________
 $(N_8,\dot{N}_5+\dot{N}_6)$ & $3$ &$b_{14}=b_{16}=...=b_{21}=0$&$b_{12}=b_{13}=b_{15}=0$\\
%___________________________________________________
 $(N_8,\dot{N}_5+\dot{N}_5)$ & $6$ &$b_{17}=b_{19}=...=b_{21}=0$&$b_{12}=...=b_{16}=b_{18}=0$\rule[-2.25mm]{0mm}{0mm}\\
%___________________________________________________
\hline
 $(N_4,\dot{N}_3+\dot{N}_5)$ & $2$ &$b_{11}=...=b_{15}=0$& $b_{9}=b_{10}=0$\rule[0mm]{0mm}{5mm}\\
%___________________________________________________
 $(N_4,\dot{N}_5+\dot{N}_5)$ & $3$&$b_{11}=b_{13}=...=0$&$b_{9}=b_{10}=b_{12}=0$\rule[-2.25mm]{0mm}{0mm}\\
\hline
 $(N_1,\dot{N}_2+\dot{N}_3)$ & $4$ &$b_{20}=...=b_{30}=0$&$b_{16}=b_{17}=b_{18}=b_{19}=0$\rule[0mm]{0mm}{5mm}\\
%_____________________________________________
 $(N_1,\dot{N}_2+\dot{N}_5)$ & $7$&$b_{23}=...=b_{30}=0$& $b_{16}=...=b_{22}=0$\\
%______________________________________________
 $(N_1,\dot{N}_2+\dot{N}_6)$ & $2$ &$b_{18}=...=b_{30}=0$& $b_{16}=b_{17}=0$\\
%___________________________________________
 $(N_1,\dot{N}_3+\dot{N}_3)$  & $7$ &$b_{22}=b_{24}=...=b_{30}=0$& $b_{16}=...=b_{21}=b_{23}=0$\\
%_____________________________________________
 $(N_1,\dot{N}_3+\dot{N}_5)$  & $10$&$b_{26}=...=b_{30}=0$&$b_{16}=...=b_{25}=0$\\
%____________________________________________
 $(N_1,\dot{N}_3+2\dot{N}_5)$  & $5$&$b_{19}=b_{22}=...=b_{30}=0$&  $b_{16}=b_{17}=b_{18}=b_{20}=b_{21}=0$\\
%______________________________________
 $(N_1,\dot{N}_3+\dot{N}_6)$  & $6$&$b_{22}=...=b_{30}=0$& $b_{16}=...=b_{21}=0$\\
%___________________________________________
 $(N_1,\dot{N}_3+\dot{N}_7)$ & $1$&$b_{17}=...=b_{30}=0$& $b_{16}=0$\\
%______________________________________________
  $(N_1,\dot{N}_4+\dot{N}_5)$ & $3$&$b_{19}=...=b_{30}=0$& $b_{16}=b_{17}=b_{18}=0$\\
%_______________________________________
 $(N_1,\dot{N}_5+2\dot{N}_3)$ & $2$&$b_{17}=b_{19}=...=b_{30}=0$& $b_{16}=b_{18}=0$\\
%_______________________________________
 $(N_1,\dot{N}_5+\dot{N}_5)$  & $11$&$b_{26}=b_{28}=...=b_{30}=0$& $b_{16}=...=b_{25}=b_{27}=0$\\
%__________________________________________
 $(N_1,3\dot{N}_5)$  & $6$& $b_{19}=b_{22}=b_{23}=b_{25}=..=b_{30}=0$&$b_{16}=..=b_{18}=b_{20}=b_{21}=b_{24}=0$\\
%___________________________________________
 $(N_1,\dot{N}_5+\dot{N}_6)$  & $8$&$b_{23}=b_{25}=...=b_{30}=0$& $b_{16}=...=b_{22}=b_{24}=0$\\ 
%_________________________________________
 $(N_1,\dot{N}_5+\dot{N}_7)$  & $4$&$b_{19}=b_{21}=...=b_{30}=0$& $b_{16}=b_{17}=b_{18}=b_{20}=0$\\
%%%%%%%%_
 $(N_1,\dot{N}_6+\dot{N}_6)$& $3$ &$b_{18}=b_{19}=b_{21}=...=b_{30}=0$& $b_{16}=b_{17}=b_{20}=0$\rule[-2.5mm]{0mm}{0mm}\\
\hline
\end{tabular}
\end{table}
\end{num}
%%%%%%%%%%%%%%%%%%%%%%%%%%%%%%%%%%%%%%%%%%%%%%%%%%%%%
\begin{num}\label{sec:delta_E_8} 
%%%%%%%%%%%%%%%%%%%%%%%%%%%%%%%%%%%%%%%%%%%%%%%%%%%%%
We end the proof with the following result:
\begin{prop}\label{prop:delta_E_8} The $\Delta$-constant strata $S_{\Delta,\calT}$ in any versal deformation of curves in $E_8$, if it is not empty, then it satisfies  
\begin{equation}\label{eq:codim_Delta}
\codim(S_{\calT},S_{\Delta,\calT})=\Delta.
\end{equation}
\end{prop}
\begin{proof} 
Since the singularities of a curve in $E_8$ out of the origin are plane curve singularities (because $E_8$ is smooth out of the origin), we can use the following classical fact for plane curves. Our proof is in the classical way of that of Proposition \ref{prop:non_trans}:
\begin{theo}\label{theo:delta_plane}[\cite{ACampo}, \cite{Tei1}, \cite{Sabir}, see \cite{Deform} Theorem 2.59] Let $T$ be a plane curve singularity with $\delta$-invariant equal to $\delta_0$. Let $F:\calD\times\CC^2\to \CC$ be its versal deformation (under $\calR$-equivalence with diffeomorphisms that move the origin). Let $\epsilon_0$ be a Milnor radius  for $F_0\inv(0)$. Take $\calD$ small enough so that $\partial B_{\epsilon_0}$ meets transversally $F_d\inv(0)$ for all $d\in \calD$. In this context, the $\delta$-constant strata is the set 
$$S_{\delta_0}:=\{d\in \calD:\sum_{p\in F_d\inv(0)\cap B_\epsilon}\delta(F_d\inv(0),p)=\delta_0\}.$$
Then:
	\begin{enumerate}
		\item the strata $S_{\delta_0}$ is contained in the closure of the following strata:%here $\delta$ singularities of $A_1$ type appear in the 0-fiber:
\begin{equation}\label{eq:A_1_curves}\{d\in\calD:F_d\inv(0)\cap B_\epsilon\ \rm{has}\ \delta_0\ \rm{singularities\  of\ type}\ A_1\}.
\end{equation} 
	\item The codimension of the strata (\ref{eq:A_1_curves}) in $\calD$ is $\delta_0$.  
	\end{enumerate}
\end{theo}

We denote by $S_{\calT;T_1,...,T_n}$ the strata of $S_{\calT}$ of curves $H_b$ with singularities of type $T_1$, ..., $T_n$ in the fiber $H_b\inv(0)\cap B_\epsilon$. Then, the stratum $S_{\Delta,\calT}$ is the union of strata $S_{\calT;T_1,...,T_n}$ with $\sum_{k=1}^n\delta(T_k)=\Delta$.

We see first that if $S_{\Delta,\calT}$ is non-empty, then 
\begin{equation}\label{eq:SS}S_{\calT;\Delta}\,\subseteq\,\overline{ S_{\calT;\underbrace{A_1,...,A_1}_{{\Delta}}}}.\end{equation} 
We prove this assertion inductively; we will see first that for any singularities $T_1$, ..., $T_n$ (with $\sum_{k=1}^n\delta(T_k)=\Delta$) we have  
\begin{equation}\label{eq:S_A_1}S_{\calT;T_1,...,T_n}\,\subseteq\, \overline{ S_{\calT;T_1,...,T_{n-1},\underbrace{A_1,...,A_1}_{{\delta(T_n)}}}}.\end{equation}
Finally, we will prove that  \begin{equation}\label{eq:codim_A_1}\codim\, S_{\calT;\underbrace{A_1,...,A_1}_\Delta}=\Delta\end{equation} 
which implies the statement of the proposition. 

To prove (\ref{eq:S_A_1}), take  a curve $\gamma$ in $\calB$ with origin in $O$ such that $\gamma(t)$ is in the stratum $S_{T_1,...,T_n}$ for any $t\neq 0$.
We can do this, using the Curve Selection Lemma. For any $t$, the curve in $X$ defined by
$H_{\gamma(t)}\inv(0)\cap B_{\epsilon}$ has precisely $n$ singularities of type $T_1$, ..., $T_n$ inside $B_\epsilon$ and outside
the origin. After a possible base change we can assume that these singularities are parametrized, say by $p_k(t)$ with $k=1,...,n$, for $p_k:(\CC,0)\to (X,O)$.

By openness of versality,
we can assume that $H|_{\calB\times B_{\epsilon}}$ is versal at any point of $\calB\times B_{\epsilon}$.

We claim that there exists a certain $N$ and a deformation of $h_0$
$$G:X\times\CC\times\CC^{N}\to\CC$$
$$(x,y,z,t,v_1,...,v_N)\mapsto G(x,y,z,t,v_1,...,v_N)$$
such that we have the equalities $$G(x,y,z,0,v_1,...,v_N)=h_0(x,y,z)\quad {\rm and}\quad    G(x,y,z,t,0,...,0)=H(x,y,z,\gamma(t))$$
for any $(v_1,...,v_N)$ and $t$ respectively, and the following additional properties:
\begin{enumerate}
\item the germ of $G_{t,v_1,...,v_N}$ at the origin has a singularity which is right equivalent
to the singularity of $G_{(t,0,...,0)}=H_{\gamma(t)}$ at the origin, for any $(v_1,...,v_N)$.
\item  the germ of $G_{t,v_1,...,v_N}$ at $p_k(t)$
has a singularity of type $T_k$ for any $k\leq n-1$, any $t\neq 0$ and any $(v_1,...,v_N)$
\item The deformation is versal at $(p_n(t),t,0,...,0)$ for any $t\neq 0$.
\end{enumerate}

The claim, together with Theorem \ref{theo:delta_plane}, implies (\ref{eq:S_A_1}). By induction we get (\ref{eq:SS}). The equality (\ref{eq:codim_A_1}) follows also by induction: using the claim and the fact that in the versal deformation of $A_1$, the $\delta$-constant stratum has codimension $1$, we get that the codimension of $S_{\calT;\underbrace{A_1,...,A_1}_k}$ in  $S_{\calT;\underbrace{A_1,...,A_1}_{k-1}}$ is $1$.

The existence of $G$ satisfying (a)-(c) is shown as in the proof of Proposition \ref{prop:non_trans}.

\end{proof}
\end{num}
%%%%%%%%%%%%%%%%%%%%%%%%%%%%%%%%%%%%%%%%%%%%%%%%%%%%%
\begin{num}
%%%%%%%%%%%%%%%%%%%%%%%%%%%%%%%%%%%%%%%%%%%%%%%%%%%%%
Taking together (\ref{eq:codim_Delta}) and (\ref{eq:codim_calA}) we get inequality (\ref{eq:ineq_codim}) and we rule out the 25 cases enumerated in Table \ref{tabla:strata} finishing the proof.
\end{num}
%%%%%%%%%%%%%%%%%%%%%%%%%%%%%%%%%%%%%%%%%%%%%%%%%%%%%
\section{Nash Problem for other quotient surface singularities.}\label{sec:quot}
%%%%%%%%%%%%%%%%%%%%%%%%%%%%%%%%%%%%%%%%%%%%%%%%%%%%
Using Theorem \ref{theo:comparacion1}, we reduce Nash Problem for quotient surface singularities to prove $\textbf{E}_8$ and $\textbf{D}_n$. The singularities $\textbf{D}_n$ can be done also using strategies of Part I and II, see \cite{tesis}. 

Anyway, we can do direct proof with our methods for most of the cases. Only $\textbf{E}_7$ -where the equations for $W^{\alpha}$ with $\alpha$ a wedge in the remaining cases have non-trivial character-, some singularities related to the dihedral $\textbf{D}_n$ in classification given in \cite{Bri} and the surface with graph as in Figure \ref{fig:non_sand} -which in particular is not a sandwiched surface singularity-, would need a more careful study. 

\begin{dibujos}
\begin{figure}[!h]
\centering
\includegraphics[width=1.6in]{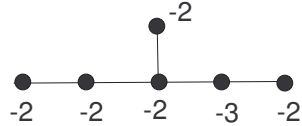}
\caption{Resolution graph of a quotient surface singularity.}\label{fig:non_sand}
\end{figure}
\end{dibujos}

Among the rest of singularities only $\textbf{E}_6$ (and $\textbf{E}_8$) needs the strategies in Part III for only one case: wedges realising the adjacency $N_1\subseteq \bar{N_3}$  with a transverse return in $N_5$ (or symmetrically, a wedge realising the adjacency $N_1\subseteq \overline{N}_5$ with a transverse return in $N_3$). See Figure \ref{fig:E_6} to see the graph. It happens that for any of such wedges, the character of the reduced equation of $W^{\underline{\alpha}}$ is also $1$; then we can proceed as in the case of $\textbf{E}_8$.

If we assume that $\textbf{E}_6$ is given by the equation $x^2+4y^3-z^5=0$ in $\CC^3$, we choose $3z^2-x$ as an equation for the image of a generic arc in $\dot{N}_1$. The stratum $S_{\calT_{3,5}}$ in the base of a versal deformation of $3z^2-x$ corresponds to the following deformation:
\begin{equation}\label{eq:versal_E_6}H_s(x,y,z)=3z^2-x+b_1x+b_2y+b_3z+b_4xy+b_5yz=0.\end{equation}
In this case, $\Delta_{1,(3,5)}=1$ and one can directly check, for instance, with \newblock {\sc Singular}, that no singularity appears in fact in $H_s\inv(0)$ for any $s$. Then the family can never be $\delta$-constant and we finish the proof.
 
\begin{dibujos}
\begin{figure}[!h]
\centering
\includegraphics[width=2.6in]{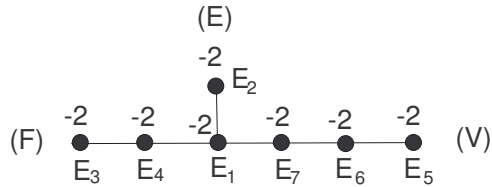}
\caption{Resolution graph of $\textbf{E}_6$.}\label{fig:E_6}
\end{figure}
\end{dibujos}

\section*{acknowledgements} I would like to thank Javier Fernandez de Bobadilla for very fruitful discussions and careful reading of this manuscript. 

I am grateful to the Algebra Department of the Faculty de Ciencias Matem\'aticas of the Universidad Complutense de Madrid for its hospitality. 

%%%%%%%%%%%%%%%%%%%%%%%%%%%%%%%%%%%
%%%%%%%%%%%%%%%%%%%%%%%%%%%%%%%%%%%%

\end{document}